\documentclass[11pt,reqno]{amsart}
\usepackage[T1]{fontenc}
\usepackage{amsfonts}
\usepackage{amssymb}
\usepackage{mathrsfs}
\usepackage[latin1]{inputenc}
\usepackage{amsmath}
\usepackage{paralist}
\usepackage[english]{babel}

\newtheorem{theorem}{Theorem}[section]
\newtheorem{lemma}[theorem]{Lemma}

\newtheorem{remark}[theorem]{Remark}
\newtheorem{prop}[theorem]{Proposition}
\newtheorem{example}[theorem]{Example}
\newtheorem{corollary}[theorem]{Corollary}

\newtheorem{hypotheses}[theorem]{Hypotheses}

\numberwithin{equation}{section}
\newcommand{\R}{{\mathbb R}}

\newcommand{\C}{{\mathbb C}}
\newcommand{\N}{{\mathbb N}}

\newcommand{\cA}{{\mathcal A}}
\newcommand{\cL}{{\mathcal L}}

\newcommand{\ve}{\varepsilon}

\newcommand{\su}{\subseteq}

\newcommand{\ov}{\overline}

\newcommand{\ar}{\mbox{\rm arg}\,}
\begin{document}

\title{Analyticity of a class of degenerate evolution equations on the canonical  simplex of $\R^d$ arising from Fleming--Viot processes}


\author{Angela A. Albanese, Elisabetta M. Mangino}



\address{ Angela A. Albanese\\
Dipartimento di Matematica ``E.De Giorgi''\\
Universit\`a del Salento- C.P.193\\
I-73100 Lecce, Italy}
\email{angela.albanese@unisalento.it}

\address{ Elisabetta M. Mangino\\
Dipartimento di Matematica ``E.De Giorgi''\\
Universit\`a del Salento- C.P.193\\
I-73100 Lecce, Italy} \email{elisabetta.mangino@unisalento.it}

\begin{abstract}
We study the analyticity of the semigroups generated by a class of  degenerate second order differential operators in the space $C(S_d)$, where $S_d$ is the canonical  simplex of $\R^d$. The semigroups  arise from the theory of Fleming--Viot processes in population genetics.
\end{abstract}

\maketitle

\section{Introduction.}

\markboth{A.\,A. Albanese, E.\, M. Mangino}%
{\MakeUppercase{Analyticity of a class of degenerate evolution equations }}

In this paper we are dealing with the class of degenerate second order elliptic differential operators 
\begin{equation}\label{e.operator}
\cA_du(x)=\frac{1}{2}\sum_{i,j=1}^dx_i(\delta_{ij}-x_j)\partial_{x_ix_j}^2u(x),\quad x\in S_d,
\end{equation}
and $m\cA_d$,
where $S_d=\{x\in [0,1]^d\mid \sum_{i=1}^dx_i\leq 1\}$ is the canonical simplex of $\R^d$
 and $m$ is a strictly positive function in the space $C(S_d)$ of all continuous functions on $S_d$.
The operator \eqref{e.operator} arises in the theory of Fleming--Viot processes as a generator of  a Markov $C_0$--semigroup defined on $C(S_d)$. Fleming--Viot processes are measure--valued processes that can be viewed as diffusion approximations of empirical processes associated with some classes of discrete time Markov chains in population genetics. We refer to \cite{EK,EK1,FV} for more details on the  topic.
In particular,  the operator \eqref{e.operator} is the generator corresponding to the diffusion model in population genetics in which neither mutation, migration, nor selection affects. This is the simplest case of a Wright--Fisher model. Actually, the generators  corresponding to more general  diffusion models in population genetics are of the following type
\begin{equation}\label{e.operatorg}
Au(x)=\frac{1}{2}\sum_{i,j=1}^dx_i(\delta_{ij}-x_j)\partial_{x_ix_j}^2u(x)+\sum_{i=1}^db_i(x)\partial_{x_i}u(x),\quad x\in S_d,
\end{equation}
where the coefficients $b_i$ belong to the space $C(S_d)$ and depend on factors as mutation, selection and migration. So, the operators (\ref{e.operatorg}) are of  degenerate  elliptic type with the elliptic part as in \eqref{e.operator}. 

 The  operators  (\ref{e.operatorg}) arising from Fleming--Viot processes have been largely studied using an analytic approach  by  several authors in different settings, see \cite{ACM-1,ACM-2,AM,CM,CR,CC,E,Met,S1,S2,S3,Stannat1} and the references quoted therein.  The interest is that the equations describing the diffusion processes are of degenerate type and hence, the classical techniques for the study of (parabolic) elliptic operators on smooth domains cannot be applied. In particular, the difficulty in studying these operators is twofold: the operators (\ref{e.operatorg}) degenerate on the boundary $\partial S_d$ of $S_d$ in a very natural way and  the domain $S_d$ is not  smooth as its boundary presents sides and corners. 
 
 As it is shown in the Feller theory for the one--dimensional case, the behaviour of the diffusion process on the boundary constitutes one of its main characteristics. So, 
 the appropriate setting for  studying  the  equations describing the diffusion process is   the space of continuous functions on the simplex $S_d$.
 
 In the one-dimensional case, the study of such type of degenerate (parabolic) elliptic problems on $C([0,1])$ started in the fifties with the papers by Feller  \cite{F1,F2}. The subsequent work of Cl\'ement and Timmermans \cite{CT} clarified which conditions on the coefficients of the operator $mA$, with $A$  defined according to (\ref{e.operatorg}) and $0<m\in C([0,1])$,  guarantee the generation of a $C_0$--semigroup in $C([0,1])$.  The problem of the regularity of the generated semigroup in $C([0,1])$ has been  considered by several authors, \cite{An,CM,BRS,Met}. In particular,  Metafune \cite{Met} established the analyticity of the semigroup under suitable conditions on the coefficients of the operator $mA$. Thus, he obtained   the analyticity of the semigroup generated by $x(1-x)D^2$ on $C([0,1])$, which was a problem left open for a long time. We refer to \cite{CMPR} for a survey on this topic.
 
In the $d$--dimensional case,  the problem of generation of a $C_0$--semigroup in $C(S_d)$ has been studied by different authors. In more generality, the problem was solved by Ethier \cite{E}. Actually. Ethier \cite{E} (see also \cite[p.375]{EK}) proved the existence of a $C_0$--semigroup of positive contractions on $C(S_d)$ under mild conditions on the drift terms $b_i$. In the following, we state such a result in the case of our interest.

\begin{theorem}\label{p.2-generation}(Ethier, \cite{E}) The closure $(\cA_d,D(\cA_d))$ of $(\cA_d,C^2(S_d))$ generates a positive and contractive $C_0$--semigroup $(T(t))_{t\geq 0}$ on $C(S_d)$. Moreover, the space $C^m(S_d)$ is a core for the infinitesimal generator of $(T(t))_{t\geq 0}$ for every $m\geq 2$. 
\end{theorem}

On the other hand, Shimakura \cite{S2} (\cite[Ch.VIII, p.221]{S3}) gave concrete representation formulas for the semigroups of diffusion processes associated to a class of Wright--Fisher models including the simplest case. In particular, Shimakura \cite[Ch.VIII, p.221]{S3} showed that 
the eigenvalues of $\cA_d$ are given by
\begin{equation}\label{e.autovalori}
\lambda_m=-\frac{m(m-1)}{2},\quad m\in\N,
\end{equation}
and that the corresponding process  is replicated on every face of $S_d$ in the following way. Denote by $\cA_{d,F}$ the restriction of $\cA_d$ to a face $F$ of $S_d$ and by $F(V)$ the face of $S_d$ having $V$ as a set of vertices. If $V$ contains $p+1$ vertices of $S_d$ with $p<d$, then $F(V)$ can be identified with the simplex $S_p$ and the differential operator $\cA_{d,F}$ with the differential operator $\cA_p$ on $S_p$, i.e.,
\begin{equation}\label{e.identi}
\cA_{d,F(V)}u=\cA_p(u|_{F(V)}),\quad u\in D(\cA_d).
\end{equation}
Moreover, in \cite[Ch.VIII, p.221]{S3} it was proved that the restriction of the semigroup  $(T(t))_{t\geq 0}$ to every  face $F(V)$ with $p+1$ vertices and $p<d$  satisfies
\begin{equation}\label{e.semig}
(T(t)f)|_{F(V)}=T_{F(V)}(t)(f|_{F(V)}),\quad f\in C(S_d), 
\end{equation}
where $(T_{F(V)}(t))_{t\geq 0}$ denotes the semigroup on $C(F(V))$ generated by $\cA_{d,F(V)}$. As the process is preserved under restriction to faces, Campiti and Rasa \cite{CR} pointed out that the domain $D(\cA_d)$ can be described recursively as follows
\begin{eqnarray}\label{e.dominio}
& & D(\cA_d)=\Big\{u\in C(S_d)\mid u\in \cap_{q\geq 1}W^{2,q}_{loc}(\stackrel{\circ}{S}_d),\ \cA_du\in C(S_d) \mbox{ and for  } \\
& & \quad  \mbox{every proper face $F\su S_d$}: u|_F\in D(\cA_{d,F}) \mbox{ and } \cA_{d,F}(u|_F)=(\cA_du)|_{F}\Big\}.\nonumber
\end{eqnarray}
If $V_d=\{v_0,\ldots, v_{d}\}$ denotes the set of vertices od $S_d$, then  \eqref{e.dominio} implies that $\cA_du(v_i)=0$ for every $u\in D(\cA_d)$ and $i=0,\ldots, d$.

There are few results about  the regularity of the generated semigroup in $C(S_d)$, \cite{ACM-1,AM}. In the papers \cite{ACM-1,AM} it was established the differentiability and the compactness of the generated semigroups related to some classes of operators of type $mA$ in $C(S_d)$,  including the generators of diffusion processes   associated to a class of  Wright--Fisher models, but not the generator \eqref{e.operator} corresponding to the simplest case. The main aim of this paper is to prove the analyticity of the semigroup generated by the closure of $(\cA_d,C^2(S_d))$ on $C(S_d)$ and hence, extending the result of Metafune to several variables.
The proof of the result is given by induction on the integer $d$. Actually, we provide a  method which  allows us to reduce the proof to the one--dimensional case and which gives information on this particular class of operators.

 The paper is organized as follows.
In \S 2 we consider the problem of the analyticity of the semigroup generated by the closure of $(\cA_2,C^2(S_2))$ on $C(S_2)$,  deepening and solving the 2-dimensional case. The end of this is to clarify in details the necessary techniques to give the inductive step.  In \S 3 we prove the analyticity of the semigroup generated by the closure of $(\cA_d,C^2(S_d))$ on $C(S_d)$ by induction.  Finally, by using the method of approximate resolvents, we show  the analyticity of the semigroup generated by the closure of $(m\cA_d,C^2(S_d))$ on $C(S_d)$.

\bigskip

\noindent{\bf Acknowledgement} The authors wish to thank Proff. A. Lunardi and G. Metafune for helpful discussions on the topic.

\bigskip

\subsection{Notation} The function spaces considered in this paper consist of complex--valued functions.

Let $K\su\R^d$ be a compact set. Denote by  $C^m(K)$ the space of all $m$--times continuously differentiable functions $u$ on $K$ such that $\lim_{x\to x_0}D^{\alpha}u(x)$ exists and is finite for all $|\alpha|\leq m$ and  $x_0\in \partial K$. In particular, $C(K)$ denotes the space of all continuous functions $u$ on $K$. The norm on $C(K)$ is the supremum norm and is  denoted by $\|\ \|_K$. The norm $\|\ \|_{m,K}$ on  $C^m(K)$  is defined by $\|u\|_{m,K}=\sum_{|\alpha|\leq m}\|D^\alpha u\|_K$.

For easy reading, in some  cases we will adopt the notation $\|\varphi(x)u\|_K$ to still denote $\sup_{x\in K}|\varphi(x)$ $u(x)|$.

A \textit{bounded analytic semigroup of angle $\theta$} with $0<\theta\leq \pi/2$ is an analytic semigroup defined in the sector $\Sigma_\theta=\{z\in \C\mid |\mbox{arg}z|<\theta\}$.


For other undefined notation and results on the theory of semigroups we refer to \cite{EN,L,P}. 

In the present paper will use some results about injective tensor products. We refer to  \cite{J,K,T,N} for definitions and  basic results  in this topic and for related applications.

\section{The $2$--dimensional case.}

\subsection{Auxiliary results}

We first consider the one--dimensional second order differential operator 
\begin{equation}\label{e.1-operator}
Au(x)=m(x)xu''(x), \quad x\in [0,b],
\end{equation}
and suppose that $b>0$ and $m$ is a strictly positive  function in $C([0,b])$. 
The operator $A$ with domain $D(A)$, defined by 
\begin{equation}\label{e.1-domainn}
D(A)=\{u\in C([0,1])\cap C^2(]0,b])\mid \lim_{x\to 0^+}Au=0, u'(b)=0\},
\end{equation}
generates a bounded analytic $C_0$--semigroup $(T(t))_{t\geq 0}$ of angle $\pi/2$ on $C([0,b])$ which is  contractive, \cite{Met,CM,CMPR,CT}.

\begin{prop}\label{p.1-gradiente} Let $b>0$ and let $m$ be a  strictly positive function in  $C([0,b])$. Then the  operator   $A$ with domain $D(A)$ defined according to  (\ref{e.1-domainn}) satisfies the  following properties.

(1) There exist $\varepsilon_b>0$, $C_b>0$ and $D_b>0$ such that,  for every $0<\varepsilon<\varepsilon_b$  and    $u\in C([0,b])\cap C^2(]0,b])$, we have 
\[ 
\|\sqrt{x} u'\|_{[0,b]} \leq \frac{C_b} {\varepsilon}  \|u\|_{[0,b]} + D_b \varepsilon \|Au\|_{[0,b]}.
 \]

(2) There exist $K_b>0$ and $t_b>0$ such that, for every $0<t<t_b$, we have
\[
\|\sqrt{x} (T(t) u)'\|_{[0,b]}\leq \frac{K_b}{\sqrt{t}}\|u\|_{[0,b]},\quad u\in C[0,b],
\]
  and such that, for every $t\geq t_b$, we have 
 \[
 \|\sqrt{x} (T(t) u)'\|_{[0,b]}\leq K_b\|u\|_{[0,b]}, \quad u\in C[0,b].
 \]
 
(3) For each $0<\theta<\pi$ there exists a costant $C_b>0$  such that, for every $\lambda\in\{z\in \C\mid |\ar z|<\theta\}$ with $|\lambda|>1$, we have
\[
\|\sqrt{x} (R(\lambda,A) u)'\|_{[0,b]}\leq \frac{C_b}{\sqrt{|\lambda|}}\|u\|_{[0,b]},\quad u\in C([0,b]).
\]
\end{prop}

\begin{proof} Denote by   $m_0=\min_{x\in [0,b]}m(x)$. Then $m_0>0$.

(1) Fix $u\in C([0,b])\cap C^2(]0,b])$. Then we have, for every $z, h\in [0,b/2]$, that 
\begin{equation}\label{eq.taylor}
 u(z+h)=u(z)+ hu'(z)+ \int_0^h(h-s)u''(z+s)ds.
 \end{equation}
Let $0<\varepsilon < \sqrt{\frac b 2}$ and $h=\varepsilon \sqrt{z}$. Then $h<\frac{b}{2}$ and  hence, from (\ref{eq.taylor}) it follows
\[ 
\sqrt{z}u'(z)=\frac 1 \varepsilon (u(z+\varepsilon\sqrt{z})-u(z)) - \frac{1}{\varepsilon}\int_0^{\varepsilon\sqrt{z}}\frac{\varepsilon \sqrt{z}-s}{z+s} u''(z+s)(z+s)ds,
\]
where
\[
 \int_0^{\varepsilon \sqrt{z}}\frac{\varepsilon\sqrt{z}-s}{z+s} \leq \frac 1 z \int_0^{\varepsilon\sqrt{z}}(\varepsilon\sqrt{z}-s)ds=\frac{\varepsilon^2}{2}.
 \]
Therefore
\begin{eqnarray*}
 \|\sqrt{z}u'\|_{[0,b/2]} &\leq &\frac{2}{\varepsilon} \|u\|_{[0,b]} + \frac{\varepsilon}{2} \|zu''\|_{[0,b]}  \leq\\
& \leq &  \frac{2}{\varepsilon}\|u\|_{[0, b]} +  \frac{\varepsilon}{2m_0} \|Au\|_{[0,b]}.
\end{eqnarray*}
On the other hand, if $z\in [b/2,b]$ and $\varepsilon \in ]0,b/4]$ (and hence, $z-\varepsilon\in [b/4,b[$), there exists $\xi\in [b/4,b]$ such that
\[
 u(z-\varepsilon)=u(z)-\varepsilon u'(z)+\frac{\varepsilon^2}{2}u''(\xi)
 \]
and hence, 
\[ 
u'(z)= \frac 1 \varepsilon (u(z)-u(z-\varepsilon)) + \frac \varepsilon 2 u''(\xi).
\]
It follows that
\begin{eqnarray*}
 |\sqrt{z}u'(z)| &\leq & \frac{2\sqrt{b}}{\varepsilon}\|u\|_{[0,b]} + \frac{\sqrt{b}}{2}\varepsilon\frac{4}{b} |\xi u''(\xi)| \\
 &\leq & \frac{2\sqrt{b}}{\varepsilon} \|u\|_{[0,b]} + \frac{2\varepsilon}{\sqrt{b}} ||zu''||_{[b/4,b]}\\
 &\leq & \frac{2\sqrt{b}}{\varepsilon} \|u\|_{[0,b]} + \frac{2\varepsilon}{\sqrt{b} m_0} \|Au\|_{[b/4,b]}.
 \end{eqnarray*}
So, 
\[ 
\|\sqrt{z}u'\|_{[b/2,b]}\leq \frac{2\sqrt{b}}{\varepsilon} \|u\|_{[0,b]} + \frac{2\varepsilon}{\sqrt{b} m_0} \|Au\|_{C[0,b]}.
\]
 We then obtain,   for every $0<\varepsilon <\varepsilon_b:=\min\{\sqrt{\frac b 2}, {\frac b 4},1\}$, that
\[ 
\|\sqrt{z} u'||_{[0,b]} \leq \frac{2+2\sqrt{b}} \varepsilon  \|u\|_{[0,b]} + \frac{\varepsilon}{m_0}(\frac 1 2 + \frac{2}{\sqrt{b}}) \|Au\|_{[0,b]}. 
\]

(2) Let $u\in C[0,b]$. Since  $(T(t))_{t\geq 0}$ is a bounded analytic $C_0$--semigroup in $C([0,b])$,  we have $T(t)u\in D(A)$ and  
there exists $M>0$ such that $t||AT(t)||\leq M$ for every $t>0$.
Applying the  property (1) above, we then obtain 
\begin{eqnarray*}
 \|\sqrt{x} (T(t)u)'\|_{[0,b]} &\leq & \frac{ C_b}{ \varepsilon}  \|T(t) u\|_{[0,1]} + D_b\varepsilon \|AT(t)u\|_{[0,b]}\\
& \leq & \frac{ C_b}{ \varepsilon}  \|u\|_{[0,b]} + D_b \varepsilon \frac{M}{t} \|u\|_{[0,b]}.
\end{eqnarray*}
Set $t_b:=\varepsilon_b^2$. Then  there exists $K_b=\max\{C_b+MD_b, C_b+\frac{MD_b}{t_b}\}>0$ such that we obtain, for every $0<t< t_b$ and taking $\varepsilon=\sqrt{t}$, that
\[ 
\|\sqrt{x} (T(t) u)'||_{[0,b]}\leq \frac{K_b}{\sqrt{t}}\|u\|_{[0,b]},
\]
 and such that,  for every  $ t\geq t_b$, 
\[ 
\|\sqrt{x} (T(t) u)'||_{[0,b]}\leq K_b \|u\|_{[0,b]}.
\]

(3) By property (2) above  the operator  $\sqrt{x}D T(t)$ is  bounded  on $C([0,b])$ with norm less or equal to $K_b/\sqrt{t}$ if  $0<t< t_b$ and to $K_b$ if $t\geq t_b$. It follows, for every $\eta>1$, $u\in C([0,b])$ and $x\in ]0,b]$, that
\[
\sqrt{x}D \left(\int_0^{+\infty} e^{-\eta t} T(t)(u)dt\right)=\int_0^\infty e^{-\eta t} (\sqrt{x} (T(t) u)'dt
 \]
and hence, 
\begin{eqnarray*}
\| \sqrt{x} (R(\eta, A)u)'\|_{[0,b]} &\leq & K_b \|u\|_{[0,b]} \left(\int_0^{t_b} t^{-1/2}e^{-\eta t}dt + \int_{t_b}^{+\infty }e^{-\eta t}dt\right) \\
&=& K_b \left(\frac{C_1}{ \sqrt{\eta}}+\frac{C_2}{\eta}\right)\|u\|_{[0,b]}
\leq \frac{C_b}{\sqrt{\eta}}\|u\|_{[0,b]}.
\end{eqnarray*}
Consequently, if $v\in D(A)$ and $\eta>1$, then 
\[
\|\sqrt{x} v'\|_{[0,b]}\leq \frac{C_b}{\sqrt{\eta}}\|\eta v- { A}v\|_{[0,b]} \leq C_b\left(\sqrt{\eta}\|v\|_{[0,b]}+\frac{1}{\sqrt{\eta}}\|{ A}v\|_{[0,b]}\right).
\]
Let $0<\theta<\pi$ be a fixed angle. If $v=R(\mu, A)u$ for some $\mu\in \{z\in \C\mid |\ar z|<\theta\}$ with $|\mu|>1$ and $u\in C([0,b])$, then by the sectoriality of  $A$ it follows
\begin{eqnarray*}
 \|\sqrt{x}(R(\mu, A)u)'\|_{[0,b]}& \leq & C_b\left(\sqrt{\eta}\|R(\mu, A)u\|_{Q_b}+\frac{1}{\sqrt{\eta}}\|AR(\mu, A)u\|_{[0,b]}\right)\\
 &=& C_b\left(\sqrt{\eta}\|R(\mu, A)u\|_{[0,b]}+\frac{1}{\sqrt{\eta}}\|\mu R(\mu, A)u-u\|_{[0,b]}\right)\\
 &\leq & C_b M\left(\frac{\sqrt{\eta}}{|\mu|}\|u\|_{[0,b]}+\frac{1}{\sqrt{\eta}}\|u\|_{[0,b]}\right),
 \end{eqnarray*}
 where the constant $M$ depends only on $\theta$.
By taking  $\eta=|\mu|$, we get the assertion.
\end{proof}

\begin{remark}\label{r.compatezza}\rm The inclusion $(D(A), \|\ \|_A)\hookrightarrow C([0,b])$ is compact and hence, $(A,D(A))$ has compact resolvent (here, $\|\ \|_A$ denotes the graph norm). Indeed, if $u\in D(A)$, then via Proposition \ref{p.1-gradiente}(1) we obtain, for every $0<x,y\leq b$, that
\begin{eqnarray}\label{e.inclu}
|u(x)-u(y)|&=&\left|\int_x^yu'(t)dt\right|=\left|\int_x^y\sqrt{t}\frac{1}{\sqrt{t}}u'(t)dt\right|\nonumber\\
&\leq & \left|\int_x^y\frac{1}{\sqrt{t}}dt\right|\|\sqrt{s}u'\|_{[0,b]}\leq C|\sqrt{x}-\sqrt{y}|\|u\|_A
\end{eqnarray}
for some  constant $C>0$. Next, let $(u_n)_n\su D(A)$ with $\sup_{n\in\N}\|u_n\|_A=K<\infty$. Then \eqref{e.inclu} implies that the sequence $(u_n)_n$ is equicontinuous in $C([0,b])$. Since $(u_n)_n$ is also equibounded in $C([0,b])$, we can apply Ascoli--Arzel\`a theorem to conclude that $(u_n)_n$ contains a subsequence
$(u_{n'})_{n'}$ converging to some $u$ in $C([0,b])$. This proves the claim.

Since $(A,D(A))$ generates an analytic $C_0$--semigroup $(T(t))_{t\geq 0}$ on $C([0,b])$ (and hence, a norm--continuous $C_0$--semigroup) and has compact resolvent, $(T(t))_{t\geq 0}$ is also compact.
\end{remark}

We now consider  the one--dimensional second order differential operator
\begin{equation}\label{e.1-model}
\cA_1u(x)=m(x)x(1-x)u''(x),\quad x\in [0,1],
\end{equation}
with domain $D(\cA_1)$ defined by
\begin{equation}\label{e.1-domains}
D(\cA_1)=\{u\in C^2(]0,1[)\cap C([0,1])\mid \lim_{x\to 0,1}\cA_1u(x)=0\}.
\end{equation}
The operator $(\cA_1,D(\cA_1))$  generates a bounded analytic $C_0$-semigroup $(T(t))_{t\geq 0}$ of angle $\pi/2$ on $C([0,1])$ which is positive and contractive, \cite{CM,Met}. Using Proposition \ref{p.1-gradiente}, we can show  that the operator $(\cA_1, D(\cA_1))$  also satisfies  the following properties.  

\begin{corollary}\label{c.1-model} Let $m$ be a  strictly positive function in  $C([0,1])$. Let $(\cA_1, D(\cA_1))$ be the differential operator on $[0,1]$ defined according to (\ref{e.1-model}).
 Then the differential operator $(\cA_1,D(\cA_1)) $  satisfies the  following properties.

(1) There exist $\ov{\varepsilon}>0$, $C>0$ and $D>0$ such that,  for every $0<\varepsilon<\ov{\varepsilon}$  and    $u\in C([0,1])\cap C^2(]0,1[)$, we have 
\[ 
\|\sqrt{x(1-x)} u'\|_{[0,1]} \leq \frac{C} {\varepsilon}  \|u\|_{[0,1]} + D \varepsilon \|\cA_1u\|_{[0,1]}.
 \]

(2) There exist $K>0$ and $\ov{t}>0$ such that, for every $0<t<\ov{t}$, we have
\[
\|\sqrt{x(1-x)} (T(t) u)'\|_{[0,1]}\leq \frac{K}{\sqrt{t}}\|u\|_{[0,1]},\quad u\in C[0,1],
\]
  and such that, for every $t\geq \ov{t}$, we have 
 \[
 \|\sqrt{x(1-x)} (T(t) u)'\|_{[0,1]}\leq K\|u\|_{[0,1]}, \quad u\in C[0,1].
 \]
 
(3) For each $0<\theta<\pi$ there exists a costant $C>0$  such that, for every $\lambda\in\{z\in\C\mid |\ar z|<\theta\}$ with $|\lambda|>1$, we have
\[
\|\sqrt{x(1-x)} (R(\lambda,\cA_1) u)'\|_{[0,1]}\leq \frac{C}{\sqrt{|\lambda|}}\|u\|_{[0,1]},\quad u\in C([0,1]).
\]
\end{corollary}

\begin{proof}(1) Let $b=\frac 1 2$. Since $m(x)(1-x)$ is a strictly positive function in $C([0,b] )$, the differential operator $\cA_1|_{[0,b]}$ is of the same type of \eqref{e.1-operator} and hence, we can apply  Proposition \ref{p.1-gradiente}(1) to conclude that there exist $\varepsilon_1>0$, $C_1>0$ and $D_1>0$ such that,  for every $0<\varepsilon<\varepsilon_1$  and    $u\in C([0,b])\cap C^2(]0,b])$, we have 
\begin{equation}\label{e.1-uno} 
\|\sqrt{x} u'\|_{[0,b]} \leq \frac{C_1}{\varepsilon}  \|u\|_{[0,b ]} + D_1 \varepsilon \|\cA_1u\|_{[0,b]}.
 \end{equation}
Next, let  $A$ be  the differential operator on $[0,b]$ defined by $Av(x)=m(1-x)x(1-x)v''(x)$ for $x\in [0,b]$,  and let 
 $\Phi\colon C([b,1])\to C([0,b])$ be the surjective isometry defined by $\Phi(u)(x):=u(1-x)$ for $u\in C([b,1])$. Then the differential operator  $A$ is of the same type of \eqref{e.1-operator}. In particular,   we have
 \[
 (A\circ \Phi)(u)(x)=m(1-x)x(1-x)u''(1-x), \quad x\in [0,b],\ u\in C([0,b]) \cap  C^2(]0,b]),
 \]
 and hence,
 \[
 (\Phi^{-1}\circ A\circ \Phi )(u)(x)=m(x)(1-x)x, \quad x\in [b,1],\ u\in C([b,1]) \cap  C^2([b,1[).
 \]
 Thus, we can apply again   Proposition \ref{p.1-gradiente}(1) to conclude that there exist $\varepsilon_2>0$, $C_2>0$ and $D_2>0$ such that,  for every $0<\varepsilon<\varepsilon_2$  and    $u\in C([b,1])\cap C^2([b,1[)$, we have  $v=\Phi (u)\in C([0,b])\cap C^2(]0,b])$ and 
 \begin{eqnarray}\label{e.1-due} 
\|\sqrt{1-x}u'\|_{[b,1]}&=&\|\sqrt{x} v'\|_{[0,b]} \leq \frac{C_2} {\varepsilon}  \|v\|_{[0,b ]} + D_2 \varepsilon \|A v\|_{[0,b]}\nonumber\\
&=&\frac{C_2}{\varepsilon}  \|u\|_{[b,1 ]} + D_2 \varepsilon \|\cA_1 u\|_{[b,1]}.
 \end{eqnarray}
Combing \eqref{e.1-uno} and \eqref{e.1-due} and setting $\ov{\varepsilon}=\min\{\varepsilon_1,\varepsilon_2\}$, we obtain, for every $0<\varepsilon<\ov{\varepsilon}$ and $u\in C([0,1])\cap C^2(]0,1[)$, that 
\begin{eqnarray*}
&  &\|\sqrt{x(1-x)}u'\|_{[0,1]} \leq \|\sqrt{x} u'\|_{[0,b]}+\|\sqrt{1-x}u'\|_{[b,1]}\\
& & \qquad \leq\frac{C_1}{\varepsilon}  \|u\|_{[0,b]} + D_1 \varepsilon \|\cA_1u\|_{[0,b]}
+\frac{C_2}{\varepsilon}  \|u\|_{[b,1]} + D_2 \varepsilon \|\cA_1 u\|_{[b,1]}\\
& &\qquad \leq\frac{C_1+C_2}{\varepsilon}  \|u\|_{[0,1]}+(D_1+D_2)\varepsilon \|\cA_1 u\|_{[0,1]}.
\end{eqnarray*}
Then, the proof of property (1) is complete.

Properties (2) and (3) follow as in the proof of Proposition \ref{p.1-gradiente}.
\end{proof}

\subsection{Consequences for a class of two--dimesional elliptic differential operators}
Using the   previous results   and  some basic properties of   injective tensor products  in the setting of Banach spaces,   \cite{J,K,T,N},  in this subsection we are able to provide resolvent estimates 
 for the two--dimensional  second order differential operators of the following type
\begin{equation}\label{e.2-operator}
A_2u(x,y)=m_1(x)x(1-x)\partial_x^2u(x,y)+ m_2(y)y\partial_y^2u(x,y),\, (x,y)\in [0,1]\times[0,b],
\end{equation}
with $b>0$, $m_1$ and $m_2$ strictly positive functions in $C([0,1])$ and in $C([0,b])$, respectively. To this end, we proceed as follows.

We  consider the one--dimensional differential operators
\[
B_1u(x)=m_1(x)x(1-x)u''(x),\ x\in [0,1], \mbox{ and }  B_2v(y)=m_2(y)yv''(y), \ y\in [0,b],
\]
with domains $D(B_1)$ and $D(B_2)$, where $D(B_1)$ is defined according to (\ref{e.1-domains}) and $D(B_2)$ is defined according   to (\ref{e.1-domainn}),
 respectively.  The operators $(B_1,D(B_1))$ and $(B_2,D(B_2))$ generate bounded analytic $C_0$--semigroups of angle $\pi/2$ on $C([0,1])$ and on $C([0,b])$ respectively, which are both  contractive. Denote  such semigroups  respectively by $(S_1(t))_{t\geq 0}$ and $(S_2(t))_{t\geq 0}$. Then the injective tensor product $(T(t))_{t\geq 0}=(S_1(t)\widehat{\otimes}_\ve S_2(t))_{t\geq 0}$ is also a bounded analytic $C_0$--semigroup of angle $\pi/2$ on $C([0,1]\times [0,b])=C([0,1])\widehat{\otimes}_\ve C([0,b])$, which is   contractive, \cite{N}. Moreover, the infinitesimal generator of $(T(t))_{t\geq 0}$  is the closure of the operator
\[
((B_1\otimes I_y) +(I_y\otimes B_2), D(B_1)\otimes D(B_2)),
\]
where $I_x$ and $I_y$ denote the identity map on $C([0,1])$ and on $C([0,b])$ with respect to the variables $x$ and $y$ respectively, and admits the space $D(B_1)\otimes D(B_2)$ as a core. Observing that
\[
A_2u=(B_1\otimes I_y)u +(I_y\otimes B_2)u,\quad u\in D(B_1)\otimes D(B_2),
\]
 we can  denote such a closure by $(A_2,D(A_2))$. Since $D(B_1)\otimes D(B_2)$ is a core for $(A_2,D(A_2))$, we have  $C^2([0,1]\times [0,b])\su D(A_2)\su C([0,1]\times [0,b])\cap C^2(]0,1[\times ]0,b])$, \cite[Chap. 44]{T}.
 
 Since the semigroups $(S_1(t))_{t\geq 0}$ and $(S_2(t))_{t\geq 0}$ are also compact, see \cite{CM,Met} and Remark \ref{r.compatezza}, their injective tensor product $(T(t))_{t\geq 0}$ shares too the compactness property, \cite[\S 44, p.285]{K}. Hence, its generator $(A_2,D(A_2))$ has compact resolvent or equivalently, the canonical injection $(D(A_2), \|\ \|_{A_2})\hookrightarrow C([0,1]\times [0,b])$ is compact, where $\|\ \|_{A_2}$ denotes the graph norm.


Next, setting $T_{2,b}=[0,1]\times [0,b]$ and using Proposition \ref{p.1-gradiente} and Corollary \ref{c.1-model}, we obtain

\begin{prop}\label{p.grest} Let $b>0$ and let $m_1$ and $m_2$ be two  strictly positive functions in $C([0,1])$ and in  $C([0,b])$ respectively. Then the  operator   $(A_2,D(A_2))$  defined according to (\ref{e.2-operator})  satisfies the  following properties.

(1) There exist $H>0$ and $\underline{t}>0$ such that, for every $0<t<\underline{t}$ and $u\in C(T_{2,b})$, we have
\[
\|\sqrt{x(1-x)} \partial_x(T(t) u)\|_{T_{2,b}}\leq \frac{H}{\sqrt{t}}\|u\|_{T_{2,b}},\ \|\sqrt{y} \partial_y(T(t) u)\|_{T_{2,b}}\leq \frac{H}{\sqrt{t}}\|u\|_{T_{2,b}} ,
\]
  and such that, for every $t\geq \underline{t}$ and $u\in C(T_{2,b})$, we have 
 \[
 \|\sqrt{x(1-x)} \partial_x(T(t) u)\|_{T_{2,b}}\leq H\|u\|_{T_{2,b}},\ \|\sqrt{y} \partial_y(T(t) u)\|_{T_{2,b}}\leq H\|u\|_{T_{2,b}}.
 \]

(2) For each $0<\theta<\pi$ there exists a costant $C>0$  such that, for every $\lambda\in\{z\in \C\mid |\ar z|<\theta\}$ with $|\lambda|>1$ and for every $u\in C(T_{2,b})$,
\[
\|\sqrt{x(1-x)}\partial_x (R(\lambda,  A_2)u)\|_{T_{2,b}} \leq  \frac{C}{\sqrt{|\lambda|}}\|u\|_{T_{2,b}},
\]
\[
\|\sqrt{y}\partial_y (R(\lambda,  A_2)u)\|_{T_{2,b}} \leq  \frac{C}{\sqrt{|\lambda|}}\|u\|_{T_{2,b}}.
\]
\end{prop}

\begin{proof} (1) By Proposition \ref{p.1-gradiente}(2) and Corollary  \ref{c.1-model} the operators  $\sqrt{x(1-x)}\partial_x S_1(t)$ and   $\sqrt{y}\partial_y S_2(t)$  are   bounded  on $C([0,1])$ and $C([0,b])$ respectively with norm less or equal to $\max\{K,K_b\}/\sqrt{t}$ if  $0<t< \underline{t}:=\min\{\ov{t},t_b\}$ and to $\max\{K,K_b\}/\sqrt{\underline{t}}$ if $t\geq \underline{t}$. Then the  operators  $(\sqrt{x(1-x)}\partial_x S_1(t)) \widehat\otimes_\varepsilon S_2(t)$ and $ S_1(t) \widehat\otimes_\varepsilon ((\sqrt{y}\partial_y S_2(t))$
are  also  bounded  on $C(T_{2,b})$ with norm less or equal to $\max\{K,K_b\}/\sqrt{t}$ if $0<t\leq \underline{t}$ or to $\max\{K,K_b\}/\sqrt{\underline{t}}$ if $t\geq \underline{t}$, \cite{J}. So, the theses follow, after having observed that, for every $u\in C(T_{2,b})$, we have
\[
\sqrt{x(1-x)} \partial_x(T(t) u)=((\sqrt{x(1-x)}\partial_x S_1(t)) \widehat\otimes_\varepsilon S_2(t))(u),
\]
\[
\sqrt{y} \partial_y(T(t) u)=(S_1(t) \widehat\otimes_\varepsilon (\sqrt{y}\partial_y S_2(t)))(u).
\]
Property (2) follows analogously to the one-dimensional case, i.e., it suffices to repeat the argument already used in the proof of Proposition \ref{p.1-gradiente}(3).
\end{proof}

We now consider  the more general case
\[
m(y)A_2u(x,y)=m(y)[m_1(x)x(1-x)\partial_x^2u(x,y)+ m_2(y)y\partial_y^2u(x,y)],\, (x,y)\in T_{2,b},
\]
and observe that

\begin{prop}\label{p.compattezza} Let $m$ be a strictly positive  function in $C([0,b])$. 
Then the operator  $(m(y) A_2, D( A_2))$ generates a contractive $C_0$--semigroup on $C(T_{2,b})$ and has compact resolvent. In particular, $D(B_1)\otimes D(B_2)$ is a core for $(m(y) A_2, D( A_2))$.
\end{prop}

\begin{proof} Since $(A_2, D(A_2))$ generates a contractive $C_0$--semigroup on $C(T_{2,b})$ and $m$ is a strictly positive  function in $C([0,b])$, we can apply a result of Dorroh \cite[Theorem]{Do} to conclude that $(m(y) A_2, D( A_2))$ also generates a contractive $C_0$--semigroup on $C(T_{2,b})$. Hence, the fact that $(m(y) A_2, D( A_2))$ has compact resolvent  follows easily, after having observed that the norms $\|\ \|_{A_2}$ and $\|\ \|_{mA_2}$ are equivalent.
\end{proof}

Thanks to Propositions \ref{p.grest} and \ref{p.compattezza} we can use the method of approximate resolvents to prove the following.

\begin{prop}\label{p.2-moperator} Let $m$ be a strictly positive  function in $C([0,b])$. 
Then the operator  $(m(y) A_2, D( A_2))$ generates an analytic   $C_0$-semigroup of angle $\pi/2$ on $C(T_{2,b})$.  The semigroup is compact.
\end{prop}

\begin{proof} For the sake of simplicity, we suppose $b=1$  and set  $m_0=:\min_{y\in [0,1]} m(y)$.  Moreover, we denote by $Q$ the square $T_{2,1}$. 

 
For each $n\in\N$ let $I_n^i:=[\frac{i-1}{n},\frac{i+1}{n}]$, $i=1, \dots,n-1$. Then we choose $\phi_n^i\in C^\infty(\R)$ for all $i=1, \dots,n-1$, such that supp\,$(\phi_n^i)\subseteq I_n^i$ and  $\sum_{i=1}^{n-1}(\phi_n^i)^2=1$. 
We observe that, if $v_i\in C(Q)$, for $i=1, \dots,n-1$,  and $y\in [0,1]$, then there exists $j\in \{1, \dots,n-1\}$ such that $y\in I_n^j$ and hence 
\[
\sum_{i=1}^{n-1}\phi_n^i(y)v_i(x,y)= \phi_n^{j-1}(y)v_{j-1}(x,y) + \phi_n^{j}(y)v_{j}(x,y)+ \phi_n^{j+1}(y)v_{j+1}(x,y).
\]
Therefore, we have 
\begin{equation}\label{parun}
\|\sum_{i=1}^{n-1}\phi_n^iv_i\|_{Q} \leq 3\sup_{i=1,\dots,n-1}\|\phi_n^iv_i\|_Q.
\end{equation}
Since the operator $(A_2,D(A_2))$ generates a bounded analytic $C_0$--semigroup of angle $\pi/2$ on $C(Q)$, for each $\lambda\in\C\setminus (-\infty,0]$, $n\in\N$ and $i=1, \dots {n-1}$, we can define
\[
 R_{in}(\lambda)=\left(\lambda - m\left(\frac{i-1}{n}\right) A_2\right)^{-1},
 \]
 and hence, for a fixed angle $0<\theta<\pi$, there exists $M>0$ such that, for every $\lambda\in \{z\in\C\mid |\ar z|<\theta \}$, $n\in\N$ and $i=1, \dots {n-1}$, we have 
\begin{equation}\label{normres}
\|R_{in}(\lambda)\|=\left[m\left(\frac{i-1}{n}\right)\right]^{-1}\left\|R\left(\frac{\lambda}{m(\frac{i-1}{n})}, A_2\right)\right\|\leq \frac{M}{|\lambda|}.
\end{equation}
%
If we set $\mu=\lambda\left[m\left(\frac{i-1}{n}\right)\right]^{-1}$, then we also have 
\begin{eqnarray*}
A_2R_{in}(\lambda)&=& \left[m\left(\frac{i-1}{n}\right)\right]^{-1}A_2R(\mu, A_2)\\
&=&\left[m\left(\frac{i-1}{n}\right)\right]^{-1}((A_2-\mu)R(\mu,A_2)+\mu R(\mu,A_2))\\
&=& \left[m\left(\frac{i-1}{n}\right)\right]^{-1}(-I+\lambda R_{in}(\lambda))
\end{eqnarray*}
and hence, 
\begin{equation}\label{normres2}
\|A_2R_{in}(\lambda)\|\leq 
 \left[m\left(\frac{i-1}{n}\right)\right]^{-1}(1+M)\leq \frac{1+M}{m_0}
\end{equation}  
We now consider the approximate resolvents of the operator $m A_2$ defined by
\[
 S_n(\lambda)u = \sum_{i=1}^{n-1} \phi_{n}^i \cdot R_{in}(\lambda) (\phi_{n}^iu), \quad u\in C(Q).
 \]
Combining (\ref{normres}) with (\ref{parun}), we obtain, for every $\lambda\in \{z\in\C\mid |\ar z|<\theta \}$  and $n\in\N$, that 
\begin{equation}\label{e.approx}
\|S_n(\lambda)\|\leq \frac{3M}{|\lambda|}.
\end{equation}
Since we have, for every $\phi,\, \eta\in D(A_2)$, that
\[ 
 A_2(\phi \eta)=\eta A_2(\phi)+ \phi A_2(\eta) + 2 [m_1(x)x\partial_x\phi\partial_x\eta + m_2(y)y \partial_y\phi\partial_y\eta],
\]  
the operators $S_n(\lambda)$ satisfy, for every $u\in C(Q)$, 
\begin{eqnarray*} 
& & (\lambda-m A_2)S_n(\lambda)u = (\lambda-m A_2)\sum_{i=1}^{n-1} \phi_n^iR_{in}(\lambda)(\phi_n^i u) \\
& &\quad =\sum_{i=1}^{n-1} \phi_n^i \cdot (\lambda-m A_2)R_{in}(\lambda) (\phi_n^iu) - \sum_{i=1}^{n-1}m A_2(\phi_n^i) \cdot R_{in}(\lambda) (\phi_n^i u) +\\
& &\quad  -2m \sum_{i=1}^{n-1} [xm_1(x)\partial_x\phi_n^i\partial_x(R_{in}(\lambda) (\phi_n^iu)) + ym_2(y) \partial_y\phi_n^i\partial_y (R_{in}(\lambda)( \phi_n^i u)]\\
& &\quad = u+ \sum_{i=1}^{n-1} \phi_n^i  \left(m\left(\frac{i-1}{n}\right)-m\right) A_2(R_{in}(\lambda) (\phi_n^iu)) +\\
& & \quad -\sum_{i=1}^{n-1} m A_2(\phi_n^i) \cdot R_{in}(\lambda) (\phi_n^i u)  -2m \sum_{i=1}^{n-1} ym_2(y) \partial_y\phi_n^i\partial_y (R_{in}(\lambda)(\phi_n^iu)) \\
& &\quad =:
(I+C_1(\lambda)+C_2(\lambda)+C_3(\lambda))u.
\end{eqnarray*}
We now fix   ${\overline{n}}\in\N$ such that $\sup_{I_{\overline{n}}^i}|m(y)-m(\frac{i-1}{{\overline{n}}})| \leq \varepsilon=:\frac{m_0}{6(1+M)}$ for $i=1,\ldots, \ov{n}-1$. Then,  from  (\ref{parun}), (\ref{normres}) and (\ref{normres2}) and Proposition \ref{p.grest}(2) it follows, for every $\lambda\in \{z\in\C\mid |\ar z|<\theta \}$ with $|\lambda|>1$ and $u\in C(Q)$,  that
\begin{eqnarray*}
& & \|C_1(\lambda)u\|_Q \leq  3\varepsilon \sup_{i=1,\dots, {\overline{n}}-1} \|\phi_{\overline{n}}^i A_2R_{i{\overline{n}}}(\lambda) (\phi_{\overline{n}}^iu)\|_Q 
 \leq  3\varepsilon\frac{1+M}{m_0} \|u\|_Q<\frac 1 2 \|u\|_Q,\\
& &\|C_2(\lambda)u\|_Q =\|\sum_{i=1}^{\ov{n}-1} m A_2(\phi_{\overline{n}}^i) \cdot R_{i{\overline{n}}}(\lambda) (\phi_{\overline{n}}^i u)\|_Q  \\
& &\leq  \sum_{i=1}^{\ov{n}-1}\max_{y\in [0,1]} m(y)\|A_2(\phi_{\overline{n}}^i)\|_Q\|R_{i{\overline{n}}}(\lambda) (\phi_{\overline{n}}^i u)\|_Q\\
& &\leq  C \max_{y\in [0,1]} m(y)  \sup_{i=1,\dots,{\overline{n}}-1}\|R_{i\overline{n}}(\lambda)(\phi_{\overline{n}}^i u)\|_Q \leq \frac{K}{|\lambda|}\|u\|_Q,\\
 & & \|C_3(\lambda)u\|_Q \leq  6 \max_{y\in [0,1]}m(y) \sup_{i=1,\dots,{\overline{n}}-1}\|ym_2(y)\partial_y\phi_{\overline{n}}^i\partial_y(R_{i\ov{n}}(\lambda)(\phi_{\ov{n}^i}u))\|_Q\\
& & \leq  H \sup_{i=1,\dots,{\overline{n}}-1}\|\sqrt{y}\partial_y(R_{i\ov{n}}(\lambda)(\phi_{\ov{n}^i}u))\|_Q\\
 & & = H\sup_{i=1,\dots,{\overline{n}}-1}\left[m\left(\frac{i-1}{{\overline{n}}}\right)\right]^{-1} \left\|\sqrt{y}\partial_y\left(R\left(\lambda m\left(\frac{i-1}{{\overline{n}}}\right)^{-1},A_2\right)(\phi_{\overline{n}}^iu)\right)\right\|_Q\\
& &\leq   \frac{K'}{\sqrt{\lambda}}\|u\|_Q,
\end{eqnarray*}
for some positive constants $K,\, K'$ independent of $\lambda$ and $u$.  
Now,  if  $|\lambda|\geq R$ for some $R>0$ large enough, then we get  $\|C_1(\lambda)+C_2(\lambda)+C_3(\lambda)\|<1$ and hence, the operator $B=(\lambda-m A_2)S_{\ov{n}}(\lambda)$ is invertible in $\cL(C(Q))$.
So, there exists $R(\lambda, m A_2)=S_{\ov{n}}(\lambda)B^{-1}$ in $\cL(Q)$  and by (\ref{e.approx}) 
\begin{equation}\label{e.rstima}
\|R(\lambda, m A_2)\|= \|S(\lambda)B^{-1}\|\leq \frac{M'}{|\lambda|}
\end{equation}
 for some  $M'>0$ independent of $\lambda$, provided $\lambda-m(y)A_2$ is injective and, in particular, for $\lambda >0$ as $m(y)A_2$ is dissipative by Proposition \ref{p.compattezza}. 
 
Observing that if $\lambda\in \rho(m(y)A_2)$ and $|\mu-\lambda|\leq \|R(\lambda, mA_2)\|^{-1}$ then $\mu\in  \rho(m(y)A_2)$, it is not difficult  to conclude via \eqref{e.rstima} and an argument of connectness that
\[
\rho(m(y)A_2)\supseteq \{z\in \C\mid |\ar z|<\theta,\ |z|>R\}.
\]
This fact togheter \eqref{e.rstima} imply that $mA_2$ generates an analytic semigroup of angle $\pi/2$.

Since the semigroup is  analytic, hence norm--continuous, and the differential operator $(mA_2, D(A_2))$  has compact resolvent,  the semigroup is also  compact. 
\end{proof}

Moreover, an analogous result of Proposition \ref{p.grest} holds in the case of the differential operator $(mA_2,D(A_2))$. Indeed, we have

\begin{prop}\label{p.2-mgradiente} Let $b>0$ and let $m$ be a strictly positive function in $C([0,b])$. Then the operator $(mA_2, D(A_2))$ satifies the following properties. 

(1) There exist $K_b>0$ and $t_b>0$ such that, for every $0<t<t_b$ and $u\in C(T_{2,b})$, we have
\[
\|\sqrt{x(1-x)}\partial_x (T(t)u)\|_{T_{2,b}} \leq \frac{K_b}{\sqrt{t}}\|u\|_{T_{2,b}},\ \ \|\sqrt{y}\partial_y (T(t)u)\|_{T_{2,b}} \leq \frac{K_b}{\sqrt{t}}\|u\|_{T_{2,b}}
\]
and such that, for every $t\geq t_b$ and $u\in C(T_{2,b})$, we have
\[
\|\sqrt{x(1-x)}\partial_x (T(t)u)\|_{T_{2,b}} \leq K_b\|u\|_{T_{2,b}}, \quad \|\sqrt{y}\partial_y (T(t)u)\|_{T_{2,b}} \leq K_b\|u\|_{T_{2,b}}.
\]

(2) For each $0<\theta<\pi$ there exist two constants $C_b>0$ and $c_b>0$   such that, for every $\lambda\in\{z\in \C\mid |\ar z|<\theta|\}$ with  $|\lambda|>c_b$ and $u\in C(T_{2,b})$, 
\[
\|\sqrt{x(1-x)}\partial_x (R(\lambda, m A_2)u)\|_{T_{2,b}} \leq  \frac{C_b}{\sqrt{|\lambda|}}\|u\|_{T_{2,b}},
\]
\[
\|\sqrt{y}\partial_y (R(\lambda, m A_2)u)\|_{T_{2,b}} \leq  \frac{C_b}{\sqrt{|\lambda|}}\|u\|_{T_{2,b}}.
\]
\end{prop}

\begin{proof} We first prove property (2). Fixed an angle $0<\theta<\pi$, let $\lambda\in \{z\in \C\mid |\ar z|<\theta|\}$ with $|\lambda|>1$ and $u\in D(m{ A}_2)=D(A_2)$. Then there exists $v\in C(T_{2,b})$ such that $R(\lambda, A_2)v=u$ and hence, 
by Proposition \ref{p.grest}(2) we have
\begin{eqnarray}\label{e.stima1}
\|\sqrt{x(1-x)}\partial_xu\|_{T_{2,b}}&=&\|\sqrt{x(1-x)}\partial_x(R(\lambda, A_2)v)\|_{T_{2,b}}\leq \frac{C}{\sqrt{|\lambda|}}\|v\|_{T_{2,b}}\nonumber\\
&=&\frac{C}{\sqrt{|\lambda|}}\|{ A}_2u-\lambda u\|_{T_{2,b}}.
\end{eqnarray}
If $|\lambda|$ is large enough, by Proposition \ref{p.2-moperator} there exists also $w\in C(T_{2,b})$ such that $R(\lambda, mA_2)w=u$ with $\|\lambda R(\lambda, mA_2)\|\leq M$ for some $M>0$. We then obtain 
\[ 
 A_2u -\lambda u = \frac 1 m (m{ A_2}u-\lambda u) + \left(\frac 1 m -1\right)\lambda u,
\] 
and hence,
 \begin{eqnarray}\label{e.stima2}
\|{ A}_2u -\lambda u\|_{T_{2,b}} &\leq & \frac{1}{ m_0} \|m{ A}_2 u-\lambda u\|_{T_{2,b}} + \left(\frac{1}{ m_0} +1\right)|\lambda| \|u\|_{T_{2,b}}\nonumber \\
&=&\frac{1}{ m_0} \|m{ A}_2 u-\lambda u\|_{T_{2,b}} + \left(\frac{1}{ m_0} +1\right)|\lambda|\|R(\lambda, mA_2)w\|_{T_{2,b}}\nonumber\\
&\leq & \frac{1}{ m_0} \|m{ A}_2 u-\lambda u\|_{T_{2,b}}+ M\left(\frac{1}{ m_0} +1\right)\|w\|_{T_{2,b}}\\
&\leq & M'\left(\frac{1}{ m_0} +1\right),\nonumber
\end{eqnarray}
with $M'=2M\left(\frac{1}{ m_0} +1\right)$ (assuming that $M\geq 1$). Combining (\ref{e.stima1}) with (\ref{e.stima2}), we 
 get 
\[
 \|\sqrt{x(1-x)}\partial_xu\|_{T_{2,b}}\leq \frac{K}{\sqrt{|\lambda|}}\|m{A}_2u-\lambda u\|_{T_{2,b}}. 
 \]
 The other inequality follows proceeding in analogous way. 
 
 (1) Fix  $0<\theta<\pi$. By  property (2) above there exist $C_b,\, c_b>0$ such that, for every $\lambda\in\{z\in \C\mid |\ar z|<\theta\}$ with $|\lambda|>c_b$ and $u\in C(T_{2,b})$,
 \[
 \|\sqrt{x(1-x)}\partial_x(R(\lambda,mA_2)u)\|_{T_{2,b}}\leq \frac{C_b}{\sqrt{|\lambda|}}\|u\|_{T_{2,b}}.
 \]
 So, we obtain, for every $\lambda\in\{z\in \C\mid |\ar z|<\theta\}$ with $|\lambda|>c_b$ and $v\in D(A_2)$, that
 \begin{eqnarray}\label{e.media}
 \|\sqrt{x(1-x)}\partial_xv\|_{T_{2,b}}& \leq &\frac{C_b}{\sqrt{|\lambda|}}\|\lambda v-mA_2v\|_{T_{2,b}}\nonumber\\
 &\leq &\frac{C_b}{\sqrt{|\lambda|}}|\lambda|\|v\|_{T_{2,b}}+\frac{C_b}{\sqrt{|\lambda|}}\|mA_2v\|_{T_{2,b}}.
 \end{eqnarray}
 Since $(mA_2,D(A_2))$ generates an analytic $C_0$--semigroup $(T(t))_{t\geq 0}$ of angle $\pi/2$ on $C(T_{2,b})$, for every $u\in C(T_{2,b})$ we have $T(t)u\in D(A_2)$ and there exist $M>0$, $w>0$  such that $t\|mA_2T(t)\|\leq Me^{wt}$ for $t>0$. Applying \eqref{e.media} with $v=T(t)u$, we then obtain that 
\begin{equation}\label{e.median}
\|\sqrt{x(1-x)}\partial_x(T(t)u)\|_{T_{2,b}}\leq C_b\sqrt{|\lambda|}\|u\|_{T_{2,b}}+\frac{C_b}{\sqrt{|\lambda|}}M\frac{e^{wt}}{t}\|u\|_{T_{2,b}}.
\end{equation}
Set $t_b:=\frac{1}{c_b}>0$. Then there exists $K_b=\max\left\{C_b(1+Me^{wt_b}), \frac{C_b}{\sqrt{t_b}}(1+M)\right\}$ such that we get, for every $0<t<t_b$ and taking $\lambda=t^{-1}$, that
\[
\|\sqrt{x(1-x)}\partial_x(T(t)u)\|_{T_{2,b}}\leq \frac{C_b}{\sqrt{t}}(1+Me^{wt_b})\|u\|_{T_{2,b}}\leq \frac{K_b}{\sqrt{t}}\|u\|_{T_{2,b}},
\]
and such that, for every $t\geq t_b$,
\[
\|\sqrt{x(1-x)}\partial_x(T(t)u)\|_{T_{2,b}}\leq \frac{C_b}{\sqrt{t_b}}(1+M)e^{wt}\|u\|_{T_{2,b}}\leq K_b e^{wt}\|u\|_{T_{2,b}}.
\]

 The other inequality follows proceeding in analogous way. 
\end{proof}

\begin{example}\label{e.model}\rm
Fix  $0<\delta <1$ and set $T_{2,{1-\delta}}=[0,1]\times [0,1-\delta]$, $T_{{1-\delta,2}}=[0,1-\delta]\times [0,1]$. Then the above results apply to  the following  second order differential operators 
\begin{eqnarray}\label{e.2-model}
A_{2,1}u(x,y)&=&\frac{1}{1-y} \frac{x(1-x)}{2} \partial_x^2u+\frac{y(1-y)}{2}\partial_y^2 u_\\
&=&\frac{1}{1-y}\left(\frac{x(1-x)}{2} \partial_x^2 u+\frac{y(1-y)^2}{2}\partial_y^2 u\right),\ (x,y)\in T_{2,{1-\delta}},\nonumber
\end{eqnarray}
with domain $D(B_1)\otimes D(B_2)$, and
\begin{eqnarray}\label{e.2-model2}
A_{2,2}u(x,y)&=&\frac{x(1-x)}{2} \partial_x^2u+\frac{1}{1-x} \frac{y(1-y)}{2}\partial_y^2 u_\\
&=&\frac{1}{1-x}\left(\frac{x(1-x)^2}{2} \partial_x^2 u+\frac{y(1-y)}{2}\partial_y^2 u\right),\quad (x,y)\in T_{{1-\delta,2}},\nonumber
\end{eqnarray}
with domain $D(B_2)\otimes D(B_1)$, where $D(B_1)$ and $D(B_2)$ are defined by 
\begin{eqnarray*}
& &D(B_1):=\{u\in C([0,1])\cap C^2(]0,1[) \mid \lim_{x\to 0^+, 1^-} x(1-x)u''(x)=0\},\\ 
& &D(B_2):=\{u\in C([0,1-\delta])\cap C^2(]0,1-\delta])\mid \lim_{y\to 0^+} yu''(y)=0, u'(1-\delta)=0\}.
\end{eqnarray*}
Indeed, as   $\delta\leq 1-y\leq 1$ for every $ y\in [0,1-\delta]$,
  by the previous considerations we can conclude that 
the closure $(A_{2,1},D(A_{2,1}))$ of  $(A_{2,1}, D(B_1)\otimes D(B_2))$ generates an analytic $C_0$-semigroup of  angle $\pi/2$ on $C(T_{2,{1-\delta}})$, which is  contractive and compact, and shares properties (1) and (2) in Proposition \ref{p.2-mgradiente}. Analogously, as $\delta\leq 1-x\leq 1$ for every $ x\in [0,1-\delta]$, the closure $(A_{2,2},D(A_{2,2}))$ of  $(A_{2,2}, D(B_2)\otimes D(B_1))$ generates an analytic $C_0$-semigroup of  angle $\pi/2$ on $C(T_{{1-\delta},2})$, which is contractive and compact, and shares properties (1) and (2) in Proposition \ref{p.2-mgradiente} with respect to $y$ and $x$. 
\end{example}

\subsection{Analyticity of a class of degenerate evolution equations on the canonical simplex of $\R^2$}
Let $S_2$ be  the simplex of $\R^2$ defined by
\[
 S_2=\{(x,y)\in\R^2\mid x,y\geq 0,\ x+y\leq 1\}.
\]
We are here concerned with the second order degenerate elliptic differential operator
\begin{equation}\label{e.2-simplex}
\cA_2u(x,y)=\frac{1}{2}x(1-x)\partial_x^2u(x,y)+\frac{1}{2}y(1-y)\partial_y^2u(x,y)-xy\partial_{xy}^2u(x,y),\quad (x,y)\in S_2.
\end{equation}


The aim of this subsection is to show the analyticity of the semigroup $(T(t))_{t\geq 0}$ generated by the closure $(\cA_2, D(\cA_2))$ of $(\cA_2, C^2(S_2))$  on $C(S_2)$ (see Theorem \ref{p.2-generation}). In order to prove this, we use suitable changes of coordinates as follows.

Fix $0<\delta<\frac 1 2$. Then,  we set
\begin{eqnarray}\label{e.2-decomposition}
\Omega_1&:=&\{(x,y)\in S_2\mid 0\leq y\leq 1-\delta\},\nonumber\\
\Omega_2&:=&\{(x,y)\in S_2\mid 0\leq x\leq 1-\delta\}.
\end{eqnarray}
Then $S_2=\cup_{i=1}^2\Omega_i$.  Next, we introduce the maps
\begin{eqnarray}\label{e.2-map}
\varphi_1\colon T_{2,{1-\delta}}\to \Omega_1, & & \quad (r,s)\to\varphi_1(r,s)=(r(1-s),s),\nonumber\\
\varphi_2\colon T_{{1-\delta},2}\to \Omega_2, & & \quad (r,s)\to\varphi_2(r,s)=(r,s(1-r)).
\end{eqnarray}

\begin{lemma}\label{l.2-continuity} The map $\varphi_i$ is bijective and a $C^\infty$--diffeomorphism for  $i=1,2$.
\end{lemma}

\begin{proof} Let $i=2$. Then, for each $(r,s)\in T_{{1-\delta},2}=[0,1-\delta]\times [0,1]$, 
$r\geq 0$, $s\geq 0$, 
and $r+s(1-r)\leq 1$ (as $(1-r)(s-1)\leq 0$ in $T_{{1-\delta},2}$), and hence, $\varphi_2(r,s)\in S_2$. Moreover, for  each $(r,s)\in T_{{1-\delta},2}$, $0\leq r\leq 1-\delta$. So, the map $\varphi_2$ is well defined. On the other hand, it is easily seen that $\varphi_2$ is continuous and injective. Next, let $(x,y)\in \Omega_2$ and let  $(r,s):=(x,\frac{y}{1-x})$. Then $(r,s)\in T_{{1-\delta},2}$ and $\varphi_2(r,s)=(x,y)$. Therefore, $\varphi_2$ is also surjective. In particular, we observe that $\varphi^{-1}(x,y)=\left(x,\frac{y}{1-x}\right)$ for $(x,y)\in \Omega_2$. Hence, $\varphi_2$ is a $C^\infty$--diffeomorphism.
The proof is analogous in the case of $\varphi_1$.
\end{proof}

So, we immediately obtain  

\begin{lemma}\label{l.2-isometry} Let $\Phi_1\colon C(\Omega_1)\to C(T_{{2,1-\delta}})$  and $\Phi_2\colon C(\Omega_2)\to C(T_{{1-\delta,2}})$ be the operators defined by
\[
\Phi_i(u)=u\circ\varphi_i,\quad u\in C(\Omega_i),\, i=1,2.
\]
Then $\Phi_i$ is a surjective isometry for  $i=1,2$. In particular, $\Phi_1(C^m(T_{2,1-\delta}))=C^m(T_{2,1-\delta})$ and $\Phi_2(C^m(T_{2,1-\delta}))= C^m(T_{1-\delta,2})$ for every $m\in\N$.
\end{lemma}

For each $i\in \{1,2\}$ we define
\begin{equation}\label{e.2-nuovi}
\cA_{2,i}:=\Phi_i^{-1}\circ A_{2,i}\circ\Phi_i, \quad D(\cA_{2,i})=\Phi_i^{-1}(D(A_{2,i})),
\end{equation}
with $(A_{2,i},D(A_{2,i}))$ the second order  differential operators defined in Example \ref{e.model}. Then the operator $(\cA_{2,i}, D(\cA_{2,i}))$ generates an analytic $C_0$--semigroup of angle $\pi/2$ on $C(\Omega_i)$ (which is also contractive and compact) for every $i\in \{1,2\}$. We observe that, for $i=1$ and $v=\Phi_1(u)$ for    some $u\in  \Phi_1^{-1}(D(A_{2,1}))$, we have
\begin{eqnarray*}
\partial_rv=(1-s)\partial_xu , \ & &	\qquad \partial_r^2v=(1-s)^2\partial_x^2u\\
\partial_sv=-r \partial_xu + \partial_yu, \ 
 & & \partial_s^2v=r^2\partial_x^2u-2r\partial_{xy}^2u+\partial_y^2u.
\end{eqnarray*}
So, we obtain 
\begin{eqnarray*}
A_{2,1}(\Phi_1(u))(r,s)&=& A_{2,1} v(r,s)= \frac{r(1-r)(1-s)}{2}\partial_x^2u + \\
& & \qquad +\frac{s(1-s)}{2}(r^2\partial_x^2u-2r\partial_{xy}^2u+\partial_{y}^2u)\\
&=& \frac{r(1-s)[1-r(1-s)]}{2}\partial_x^2u-rs(1-s)\partial_{xy}^2u+\\
& & +\qquad \frac{s(1-s)}{2}\partial_y^2u,\qquad (r,s)\in T_{2,{1-\delta}},
\end{eqnarray*}
and hence, 
\[
\cA_{2,1}(x,y)=\frac{x(1-x)}{2}\partial_x^2u(x,y)+\frac{y(1-y)}{2}\partial_y^2u(x,y)-xy\partial_{xy}u(x,y),\ (x,y)\in \Omega_1,
\]
i.e., $\cA_2|_{\Omega_1}=\cA_{2,1}$. 

On the other hand, for $i=2$ and $v=\Phi_2(u)$ for   some $u\in  \Phi_2^{-1}(D(A_{2,2}))$, we have
\begin{eqnarray*}
\partial_rv=\partial_xu -s\partial_yu, \ & &	\qquad \partial_r^2v=\partial_x^2u-2s\partial_{xy}^2u+s^2\partial_y^2u\\
\partial_sv=(1-r) \partial_yu , \ 
 & & \partial_s^2v=(1-r)^2\partial_y^2u.
\end{eqnarray*}
So, we obtain 
\begin{eqnarray*}
A_{2,2}(\Phi_2(u))(r,s)&=& A_{2,2} v(r,s)= \frac{r(1-r)}{2}(\partial_x^2u -2s\partial^2_{xy}u+s^2\partial_{y}^2u)+ \\
& & \qquad + \frac{s(1-s)(1-r)}{2}\partial_{y}^2u)\\
&=& \frac{r(1-r)}{2}\partial_x^2u-rs(1-s)\partial_{xy}^2u+\\
& & \quad +\frac{s(1-r)[1-s(1-r)]}{2}\partial_y^2u,\qquad (r,s)\in T_{{1-\delta},2},
\end{eqnarray*}
and hence,  
\[
\cA_{2,2}(x,y)=\frac{x(1-x)}{2}\partial_x^2u(x,y)+\frac{y(1-y)}{2}\partial_y^2u(x,y)-xy\partial_{xy}u(x,y),\ (x,y)\in \Omega_2,
\]
i.e., $\cA_2|_{\Omega_2}=\cA_{2,2}$. 

We may now prove the main theorem of this section.

\begin{theorem}\label{t.2-analyticity}
The closure $(\cA_2, D(\cA_2))$ of $(\cA_2, C^2(S_2))$ generates an  analytic  $C_0$--semigroup $(T(t))_{t\geq 0}$ of angle $\pi/2$ on $C(S_2)$. The semigroup is compact.
\end{theorem}

\begin{proof}
Fix $0<\delta<\frac 1 2$. 
Let $\{\psi_i\}_{i=1,2}\su C_c^\infty(\R^2)$ such that $\sum_{i=1}^2(\psi_i)^2=1$ on $S_2$ and 
\begin{eqnarray*}
& & \mathrm{supp}(\psi_1) \subseteq \{ (x,y)\in \R^2\mid y< 1-\delta\}\\
& &\mathrm{supp}(\psi_2)\subseteq\{ (x,y)\in \R^2\mid x<1-\delta\}.
\end{eqnarray*}
For the sake of simplicity, 
we still denote by $\psi_i$ the restriction of $\psi_i$ to $\Omega_i$, for $i=1,2$. 

By Proposition \ref{p.2-moperator} the operators $\cA_{2,i}$ generate analytic $C_0$--semigroups of angle $\pi/2$. So, if $0<\theta<\pi$ is a fixed angle, we can find two positive constants $C$ and $R$ such that, for $|\lambda|\geq R$ with $|\ar \lambda|<\theta$, the resolvents $R(\lambda, \cA_{2,i})$ exist and satisfy
\begin{equation}\label{e.risolvente}
\|R(\lambda, \cA_{2,i})\|\leq \frac{C}{|\lambda|},\quad i=1,2.
\end{equation} 
So, 
for every $\lambda\in\{z\in\C\mid |\ar z|<\theta\}$ with $|\lambda|\geq R$, we  can define the operator $S(\lambda)\colon C(S_2)\to C(S_2)$ via  
\begin{equation}\label{e.2-resolventt}
S(\lambda)u= \sum_{i=1}^2 \psi_i R(\lambda, \cA_{2,i}) (\psi_i u),\quad u\in C(S_2),
\end{equation}
and hence, 
\[
\|S(\lambda)\|\leq \frac{3C}{|\lambda|}.
\]
We observe that  the previous considerations on the differential operators $\cA_{2,i}$ ensure,  for every $i=1,2$ and $u\in C(S_2)$,  that
\begin{equation}\label{e.2-compr}
 \cA_2(\psi_i R(\lambda, \cA_{2,i})(\psi_iu))= \cA_{2,i}(\psi_iR(\lambda, \cA_{2,i})(\psi_i u)), 
 \end{equation}
and, for every $f,g\in D(\cA_{2,i})$, that 
\begin{eqnarray}\label{e.2-compp} 
\cA_{2,i}(fg)&=&(\cA_{2,i}f)g+ f(\cA_{2,i}g)+ [x(1-x-y)\partial_xf\partial_xg + \\
& &  +(1-x-y)y \partial_yf\partial_yg+xy(\partial_xf-\partial_yf)(\partial_xg-\partial_yg)]\nonumber
\end{eqnarray}
By (\ref{e.2-compr}) and (\ref{e.2-compp}) we obtain, for every  $\lambda\in \{z\in\C\mid |\ar z|<\theta\}$ with $|\lambda|\geq R$ and $u\in C(S_2)$, that
\begin{eqnarray*}
(\lambda-\cA_2)S(\lambda)u &=& \lambda S(\lambda)(u)- \sum_{i=1}^2 \cA_2(\psi_i R(\lambda, \cA_{2,i})(\psi_iu)) \\
&=&\lambda S(\lambda)(u)-\sum_{i=1}^2 \cA_{2,i}(\psi_i(R(\lambda, \cA_{2,i})(\psi_i u))\\
&=&\sum_{i=1}^2\psi_i (\lambda-\cA_{2,i})R(\lambda, \cA_{2,i})(\psi_i u)- \sum_{i=1}^2 \cA_{2,i}(\psi_i) R(\lambda, \cA_{2,i})(\psi_i u) +\\
&-& \sum_{i=1}^2 [x(1-x-y)\partial_x(R(\lambda,\cA_{2,i})(\psi_iu))\partial_x\psi_i+\\
&+& y(1-x-y)\partial_y(R(\lambda,\cA_{2,i})(\psi_iu))\partial_y\psi_i+\\
&-& xy (\partial_x(R(\lambda,\cA_{2,i})(\psi_iu))- \partial_y(R(\lambda,\cA_{2,i})(\psi_iu)))(\partial_x\psi_i-\partial_y\psi_i)]\\
&=:& (I + B(\lambda) + C(\lambda))(u).
\end{eqnarray*}
Applying \eqref{e.risolvente} we obtain, for every $\lambda\in\{z\in\C\mid |\ar z|<\theta\}$ with $ |\lambda|\geq R$ and $u\in C(S_2)$, that 
 \[
 \|B(\lambda) u\|_{S_2}\leq \frac{M}{|\lambda|} \|u\|_{S_2},
 \]
 for some $M>0$.
In order to estimate $C(\lambda)$ we proceed as  follows.

Let $f\in D(\cA_{2,1})$ and  $v=\Phi_1(f)$. Then $v\in D(A_{2,1})$ and the following  holds   
\begin{equation}\label{e.2-derivate3}
 \partial_rv= (1-s)\partial_xf,\quad \partial_sv=-r \partial_xf + \partial_yf,
 \end{equation}
 and hence,
\begin{equation}\label{e.2-cderivate3}
 \partial_xf= \frac{\partial_rv}{1-s},\ \partial_yf=\partial_sv+\frac{r}{1-s}\partial_rv,\ \partial_xf-\partial_yf=\frac{1-r}{1-s}\partial_rv -\partial_sv.
 \end{equation}
So, by Lemma \ref{l.2-isometry} and Proposition \ref{p.2-mgradiente}(2) (combined  with Example \ref{e.model}) we obtain, for every $\lambda\in \{z\in\C\mid |\ar z|<\theta\}$ with $|\lambda|>R_1:=\max\{R,c_b\}$, that     
\begin{eqnarray*}
& &\|x(1-x-y)\partial_xf\partial_x\psi_1+ y(1-x-y) \partial_yf\partial_y\psi_1+ xy (\partial_xf-\partial_yf)(\partial_x\psi_1-\partial_y\psi_1)\|_{S_2}\\
& &\leq \|x(1-x-y)\partial_xf\partial_x\psi_1\|_{\Omega_1}+\|y(1-x-y) \partial_yf\partial_y\psi_1\|_{\Omega_1}+\\
& &\qquad +\|xy (\partial_xf-\partial_yf)(\partial_x\psi_1-\partial_y\psi_1)\|_{\Omega_1}\\
& &\leq C_1\Big(\left\|r(1-s)(1-r)(1-s)\frac{v_r}{1-s}\right\|_{T_{2,{1-\delta}}}+\left\|(1-r)(1-s)s\left(v_s + \frac{r}{1-s}v_r\right)\right\|_{T_{2,{1-\delta}}} \\
& &+ \left\|r(1-s)s \left(v_r\frac{1-r}{1-s}-v_s\right)\right\|_{T_{2,1-\delta}}\Big)\\
& & \leq\frac{C'_1}{\sqrt{|\lambda|}} \|\lambda v-A_{2,1} v\|_{T_{2,{1-\delta}}}= \frac{C'_1}{\sqrt{|\lambda|}}\|\lambda f -\cA_{2,1} f\|_{\Omega_1}.
\end{eqnarray*}
If $f=R(\lambda, \cA_{2,1})(\psi_1u)$ for some $u\in C(S_2)$, then  it follows, for every $\lambda\in\{z\in\C\mid |\ar z|<\theta\}$ with $|\lambda|>R_1$, that 
\begin{eqnarray}\label{e.2-rsol3}
 & & \|x(1-x-y)\partial_x(R(\lambda,\cA_{2,1})(\psi_1u))\partial_x\psi_1+ y(1-x-y)\partial_y(R(\lambda,\cA_{2,1})(\psi_1u))\partial_y\psi_1+\nonumber\\
& &+ xy (\partial_x(R(\lambda,\cA_{2,1})(\psi_1u))- \partial_y(R(\lambda,\cA_{2,1})(\psi_1u)))(\partial_x\psi_1-\partial_y\psi_1)\|_{S_2}\nonumber\\
& & \quad \leq \frac{C_1'}{\sqrt{|\lambda|}}\|u\|_{\Omega_1}\leq \frac{C_1'}{\sqrt{|\lambda|}}\|u\|_{S_2},
\end{eqnarray}
with $C_1'$ a positive constant independent of $\lambda$ and $f$. 

By the simmetry of the change of variables, one analogously shows that, there exists $C_2>0$ such that, for every $\lambda\in\{z\in\C\mid |\ar z|<\theta\}$ with $|\lambda|>R_1$ and $u\in C(S_2)$, we have
\begin{eqnarray}\label{e.2-rsol2}
 & & \|x(1-x-y)\partial_x(R(\lambda,\cA_{2,2})(\psi_2u))\partial_x\psi_2+ (1-x-y)y \partial_y(R(\lambda,\cA_{2,2})(\psi_2u))\partial_y\psi_2\nonumber\\
& & +xy (\partial_x(R(\lambda,\cA_{2,2})(\psi_2u))- \partial_y(R(\lambda,\cA_{2,1})(\psi_2u)))(\partial_x\psi_2-\partial_y\psi_2)\|_{S_2}\nonumber\\
& & \quad \leq  \frac{C_2}{\sqrt{|\lambda|}}\|u\|_{S_2}.
\end{eqnarray}

Combining  (\ref{e.2-rsol2}) and (\ref{e.2-rsol3}), we obtain that there exists $K>0$ such that, for every $\lambda\in\{z\in\C\mid |\ar z|<\theta\}$ with $|\lambda|>R_1$ and $u\in C(S_2)$,
\[
\|C(\lambda)u\|_{S_2} \leq \frac{K}{\sqrt{|\lambda|}}\|u\|_{S_2}.
\]
If $|\lambda|>>1$, then the operator $B=(\lambda-\cA_2)S(\lambda)$ is invertible in $\cL(C(S_2))$. So, there exists $R(\lambda, \cA_2)=S(\lambda)B^{-1}$ and
\[
\|R(\lambda, \cA_2)\|=\|S(\lambda)B^{-1}\| \leq \frac{C'}{|\lambda|},
\]
with $C'$ a positive constant independent of $\lambda$, provided $\lambda-\cA_2$ is injective and, in particular, for $\lambda>0$ as $\cA_2$ is dissipative. To conclude that the semigroup is analytic of angle $\pi/2$ it now suffices to repeat the argument already used in  the proof of Proposition \ref{p.2-moperator}.

Since $R(\lambda, \cA_2)=S(\lambda)B^{-1}$ for some $\lambda>0$ and the operator $S(\lambda)$ is compact by Proposition \ref{p.2-moperator}, the differential operator $(\cA_2,D(\cA_2))$ has compact resolvent. Thus, the semigroup is also compact, being analytic and hence, norm continuous.
\end{proof}

Recalling that the eigenvalues of the operator $\cA_2$ are given by $\lambda_m=-\frac{m(m-1)}{2}$, $m\geq 1$, \cite[Ch.VIII, p.221]{S3}, and using Theorem \ref{t.2-analyticity} together with \cite[Proposition 5.6]{Met} we obtain the following result.

\begin{theorem}\label{t.2-boundedanal} The semigroup generated by $(\cA_2, D(\cA_2))$ is bounded analytic of angle $\pi/2$.
\end{theorem}

Moreover, the differential operator $(\cA_2, D(\cA_2))$ satisfies 

\begin{prop}\label{p.2-gradiente} The closure $(\cA_2, D(\cA_2))$  of the differential operator $(\cA_2, C^2(S_2))$ defined in (\ref{e.2-simplex}) satisfies the following properties.

(1) There exist $K_b>0$ and $t_b>0$ such that, for every $0<t<t_b$ and $u\in C(S_2)$, we have
\[
\|\sqrt{x(1-x)}\partial_x (T(t)u)\|_{S_{2}} \leq \frac{K_b}{\sqrt{t}}\|u\|_{S_{2}},\ \ \|\sqrt{y(1-y)}\partial_y (T(t)u)\|_{S_{2}} \leq \frac{K_b}{\sqrt{t}}\|u\|_{S_{2}}
\]
and such that, for every $t\geq t_b$ and $u\in C(S_{2})$, we have
\[
\|\sqrt{x(1-x)}\partial_x (T(t)u)\|_{S_{2}} \leq K_b\|u\|_{S_{2}}, \quad \|\sqrt{y(1-y)}\partial_y (T(t)u)\|_{S_{2}} \leq K_b\|u\|_{S_{2}}.
\]

(2) For each $0<\theta<\pi$ there exists two constants $C>0$ and $l>1$ such that, for every $\lambda\in \{z\in \C\mid |\ar z|<\theta \}$ with $|\lambda|>l$ and $u\in C(S_2)$, we have
\[
\|\sqrt{x(1-x)}\partial_x(R(\lambda, \cA_2)u)\|_{S_2}\leq \frac{C}{\sqrt{|\lambda|}}\|u\|_{S_2},
\]
\[
\|\sqrt{y(1-y)}\partial_y(R(\lambda, \cA_2)u)\|_{S_2}\leq \frac{C}{\sqrt{|\lambda|}}\|u\|_{S_2}.
\]
\end{prop}

\begin{proof} We first prove property (2). According to the notation in the proof of Theorem \ref{t.2-analyticity}, fixed an angle $0<\theta<\pi$  there exists $l>1$ such that,  for every $\lambda\in \{z\in\C\mid |\ar z|<\theta\}$ with $|\lambda|>l$, we have
\[
R(\lambda, \cA_2)=S(\lambda)B^{-1},
\]
where the operators $S(\lambda)$ are defined according to (\ref{e.2-resolventt})  and $\|B^{-1}\|\leq 2$.
So, we obtain,  for every $\lambda\in \{z\in\C\mid |\ar z|<\theta\}$ with $|\lambda|>l$ and $u\in C(S_2)$,  that
\begin{eqnarray}\label{e.2-rgradiente}
& &\|\sqrt{x(1-x)}\partial_x(R(\lambda, \cA_2)u)\|_{S_2} \leq  \sum_{i=1}^2 \|\sqrt{x(1-x)}\partial_x(\psi_i R(\lambda, \cA_{2,i}) (\psi_i B^{-1}u))\|_{\Omega_i}\nonumber\\
& & \leq  \sum_{i=1}^2 \Big(\|\sqrt{x(1-x)}\partial_x(\psi_i) R(\lambda, \cA_{2,i}) (\psi_i B^{-1}u)\|_{\Omega_i}+\nonumber\\
& & \qquad +\|\sqrt{x(1-x)}\psi_i \partial_x(R(\lambda, \cA_{2,i}) (\psi_i B^{-1}u))\|_{\Omega_i}\Big)\\
& & \leq  c \sum_{i=1}^2\left(\frac{C}{|\lambda|}\|\psi_i B^{-1}u\|_{\Omega_i}+\|\sqrt{x(1-x)}\psi_i \partial_x(R(\lambda, \cA_{2,i}) (\psi_i B^{-1}u))\|_{\Omega_i}\right)\nonumber\\
&  & \leq  c \sum_{i=1}^2\left(\frac{C}{\sqrt{|\lambda|}}\|B^{-1}\|\|u\|_{S_2}+\|\sqrt{x(1-x)} \partial_x(R(\lambda, \cA_{2,i}) (\psi_i B^{-1}u))\|_{\Omega_i}\right),\nonumber
\end{eqnarray}
where $c:=\sup_{i=1}^2\|\psi_i\|_{1,\Omega_i}$. To estimate the second addend on the right in (\ref{e.2-rgradiente}) we proceed as follows.

For $i=1$ set $f= R(\lambda, \cA_{2,1}) (\psi_1 B^{-1}u)$ and $v=\Phi_1(f)$. Then by (\ref{e.2-derivate3}) and (\ref{e.2-cderivate3})
\begin{eqnarray}\label{e.2-adde3}
\|\sqrt{x(1-x)}\partial_xf\|_{\Omega_1}&=&\left\|\frac{\sqrt{r(1-s)[1-r(1-s)]}}{1-s}\partial_rv\right\|_{T_{2,1-\delta}}\nonumber\\
&\leq &\frac{C''_1}{\sqrt{|\lambda|}}\|\lambda v-A_{2,1}v\|_{T_{2,{1-\delta}}}=\frac{C_1}{\sqrt{|\lambda|}}\|\lambda f-\cA_{2,1}f\|_{\Omega_1}\nonumber\\
&=&\frac{C_1}{\sqrt{|\lambda|}}\|\psi_1 B^{-1}u\|_{\Omega_1}\leq \frac{C_1}{\sqrt{|\lambda|}}\|B^{-1}\|\|u\|_{S_2}.
\end{eqnarray} 
Next, for $i=2$ set $f= R(\lambda, \cA_{2,2}) (\psi_2 B^{-1}u)$ and $v=\Phi_2(f)$. Then $\partial_xf=\partial_rv+\frac{s}{1-r}\partial_sv$ and hence,
\begin{eqnarray}\label{e.2-adde3'}
\|\sqrt{x(1-x)}\partial_xf\|_{\Omega_2}&=&\left\|\sqrt{r(1-r)}\left(\partial_rv+\frac{s}{1-r}\partial_sv\right)\right\|_{T_{1-\delta,2}}\nonumber\\
&\leq & \left\|\sqrt{r(1-r)}\partial_rv\right\|_{T_{1-\delta,2}}+\left\|\sqrt{r(1-r)}\frac{s}{1-r}\partial_sv\right\|_{T_{1-\delta,2}}\nonumber\\
&\leq &\frac{C''_1}{\sqrt{|\lambda|}}\|\lambda v-A_{2,2}v\|_{T_{2,{1-\delta}}}=\frac{C_2}{\sqrt{|\lambda|}}\|\lambda f-\cA_{2,2}f\|_{\Omega_2}\nonumber\\
&=&\frac{C_2}{\sqrt{|\lambda|}}\|\psi_1 B^{-1}u\|_{\Omega_1}\leq \frac{C_2}{\sqrt{|\lambda|}}\|B^{-1}\|\|u\|_{S_2}.
\end{eqnarray} 
Combining \eqref{e.2-rgradiente}, \eqref{e.2-adde3} and \eqref{e.2-adde3'}  we obtain, for every  $\lambda\in \{z\in\C\mid |\ar z|<\theta\}$ with $|\lambda|>l$ and $u\in C(S_2)$,  that
\[
\|\sqrt{x(1-x)}\partial_x(R(\lambda, \cA_2)u)\|_{S_2} \leq \frac{C}{\sqrt{|\lambda|}}\|u\|_{S_2}
\]
for some constant $C>0$ independent of $u$ and $\lambda$.

The proof of the other case is analogue.

Property (1) follows  as in the proof of Proposition \ref{p.2-mgradiente}(1). 
\end{proof}

\section{The $d$--dimensional case.}
The aim of this section is to show that the semigroup $(T(t))_{t\geq 0}$ generated by the closure  $(\cA_d, D(\cA_d))$ of the operator $(\cA_d, C^2(S_d))$ on $C(S_d)$ (see Theorem \ref{p.2-generation}) is also analytic for $d>2$. We  prove this using an argument by induction as follows.

\subsection{Inductive Hypotheses and Consequences}


 We suppose that the following holds.

\begin{hypotheses}[\bf Inductive Hypothesis]\label{Hypo} Suppose that the closure $(\cA_d, D(\cA_d))$ of $(\cA_d,C^2(S_d))$ satisfies the following properties.

(1) $(\cA_d, D(\cA_d))$ generates a bounded analytic $C_0$--semigroup $(T(t))_{t\geq 0}$ of angle $\pi/2$ on $C(S_d)$. The semigroup is compact.

(2) For each $0<\theta<\pi$ there exist  $C>0$ and $l>1$ such that, for every $\lambda\in\{z\in \C\mid |\ar z|<\theta\}$ with $|\lambda|>l$, $i=1,\ldots, d$  and $u\in C(S_d)$, we have
\[
\|\sqrt{x_i(1-x_i)}\partial_{x_i}(R(\lambda, \cA_{d})u)\|_{S_d}\leq \frac{C}{\sqrt{|\lambda|}}\|u\|_{S_d}.
\]
\end{hypotheses}

 In order to prove the inductive step, we need to provide some auxiliary results as follows.

Fix $0<\delta<\frac 1 2$. Then, we define  the sets
$T_{d+1,1-\delta}:=S_d\times [0,1-\delta]$,  $T_{1-\delta, d+1}:=[0,1-\delta]\times S_{d}$,
and consider the second order differential operators
\begin{equation}\label{e.d-operatore1}
A_{d+1,1}=\frac{1}{2(1-x_{d+1})}\sum_{i,j=1}^dx_i(\delta_{ij}-x_j)\partial_{x_ix_j}^2u+\frac{1}{2}x_{d+1}(1-x_{d+1})\partial_{x_{d+1}}^2u, \ x\in T_{d+1,1-\delta}
\end{equation}
\begin{equation}\label{e.d-operatore2}
A_{d+1,2}=\frac{1}{2}x_{1}(1-x_{1})\partial_{x_{1}}^2u+\frac{1}{2(1-x_{1})}\sum_{i,j=2}^{d+1}x_i(\delta_{ij}-x_j)\partial_{x_ix_j}^2u, \ x\in T_{1-\delta,d+1}.
\end{equation}
Recalling that $C(T_{d+1,1-\delta})=C(S_d)\widehat{\otimes}_\ve C([0,1-\delta])$ and $C(T_{1-\delta,d+1})=C([0,1-\delta])\widehat{\otimes}_\ve C(S_d)$ and the discussion prior to Proposition \ref{p.grest}, we can prove via  Hypotheses \ref{Hypo}, Proposition \ref{p.1-gradiente} and analogously to the proofs of Propositions \ref{p.grest}, \ref{p.2-moperator} and \ref{p.2-mgradiente} the following facts.

Setting $D(B):=\{u\in C([0,1-\delta])\cap C^2(]0,1-\delta])\mid \lim_{y\to 0^+} yu''(y)=0,\ u'(1-\delta)=0\}$, we have
 
\begin{prop}\label{p.d-anal} Suppose that Hypotheses \ref{Hypo} hold. Then the following properties are satisfied.

(1) The closure $(A_{d+1,1},D(A_{d+1,1}))$ of $(A_{d+1,1}, D(\cA_d)\otimes D(B))$ generates an analytic $C_0$--semigroup of angle $\pi/2$ on $C(T_{d+1,1-\delta})$ which is  contractive and compact.

(2) The closure $(A_{d+1,2},D(A_{d+1,2}))$ of $(A_{d+1,2},  D(B)\otimes D(\cA_d))$ generates an analytic $C_0$--semigroup of angle $\pi/2$ on $C(T_{1-\delta,d+1})$ which is contractive and compact.
\end{prop}

\begin{prop}\label{p.d-stimerisol} Suppose that Hypotheses \ref{Hypo} hold. Then the following properties are satisfied.

(1) For each $0<\theta<\pi$ there exist  $C_1>0$ and $l_1>1$ such that, for every $\lambda\in\{z\in \C\mid |\ar z|<\theta\}$ with $|\lambda|>l_1$ and $u\in C(T_{d+1,1-\delta})$, we have
\[
\|\sqrt{x_i(1-x_i)}\partial_{x_i}(R(\lambda, A_{d+1,1})u)\|_{T_{d+1,1-\delta}}\leq \frac{C_1}{\sqrt{|\lambda|}}\|u\|_{T_{d+1,1-\delta}},\quad i=1,\ldots, d,
\]
and
\[
\|\sqrt{x_{d+1}}\partial_{x_{d+1}}(R(\lambda, A_{d+1,1})u)\|_{T_{d+1,1-\delta}}\frac{C_1}{\sqrt{|\lambda|}}\|u\|_{T_{d+1,1-\delta}}.
\]

1) For each $0<\theta<\pi$ there exist  $C_2>0$ and $l_2>1$ such that, for every $\lambda\in\{z\in \C\mid |\ar z|<\theta\}$ with $|\lambda|>l_2$ and $u\in C(T_{1-\delta,d+1})$, we have
\[
\|\sqrt{x_{1}}\partial_{x_{1}}(R(\lambda, A_{d+1,2})u)\|_{T_{1-\delta,d+1}}\frac{C_2}{\sqrt{|\lambda|}}\|u\|_{T_{1-\delta,d+1}},
\]
and
\[
\|\sqrt{x_i(1-x_i)}\partial_{x_i}(R(\lambda, A_{d+1,2})u)\|_{T_{1-\delta,d+1}}\leq \frac{C_2}{\sqrt{|\lambda|}}\|u\|_{T_{1-\delta,d+1}}, \quad i=2,\ldots, d+1.
\]
\end{prop}

\subsection{Analyticity of a class of degenerate evolution equations on the canonical simplex of $\R^{d+1}$}
For the  inductive step, we also need to 
  perfom  the following  changes of coordinates.
 
 Fix $0<\delta<\frac 1 2$. Then,  we set
\begin{eqnarray}\label{e.d-decomposition}
& &\Omega_{1}:=\left\{x\in S_{d+1}\mid 0\leq x_{d+1}\leq 1-\delta\right\},\nonumber\\
& &\Omega_2:=\left\{x\in S_{d+1}\mid 0\leq x_1\leq 1-\delta\right\}.
\end{eqnarray}
Then $S_{d+1}=\cup_{i=1}^{2}\Omega_i$.  Next, we  consider the maps $\varphi_1\colon T_{d+1,1-\delta}\to \Omega_1$ and $\varphi_2\colon T_{1-\delta,d+1}\to \Omega_2$ defined by
\begin{eqnarray}\label{e.d-map}
& & \varphi_{1}(r):=(r_1(1-r_{d+1}), r_2(1-r_{d+1}),\ldots,r_d(1-r_{d+1}),r_{d+1})\nonumber\\
& &\varphi_{2}(r):=(r_1, r_2(1-r_{1}),\ldots,r_d(1-r_{1}),r_{d+1}(1-r_{1}).
\end{eqnarray}

\begin{lemma}\label{l.d-continuity}  The map $\varphi_i$ is  bijective and a $C^\infty$--diffeomorphism for $i=1,2$.
\end{lemma}

\begin{proof} 
The map $\varphi_{1}$ is well defined. Indeed, for each $r\in T_{d+1,1-\delta}$, $r_i\geq 0$ for $i=1,\ldots,d+1$ and $\sum_{i=1}^dr_i(1-r_{d+1})+r_{d+1}\leq 1$ (as $(1-r_{d+1})(\sum_{i=1}^d-1)\leq 0$ in $T_{d+1,1-\delta}$) and hence, $\varphi_1(r)\in S_{d+1}$. Moreover, for each $r\in T_{d+1,1-\delta}$, $0\leq r_{d+1}\leq 1-\delta$. So, the map $\varphi_1$ is well defined. On the other hand,  $\varphi_1$ is clearly continuous and injective. Next, let $x\in \Omega_1$ and set $r=\left(\frac{x_1}{1-x_{d+1}}, \frac{x_2}{1-x_{d+1}},\ldots, \frac{x_d}{1-x_{d+1}}, x_{d+1}\right)$. Then $x\in T_{d+1,1-\delta}$ and $\varphi_1(r)=x$. Thus, $\varphi_1$ is also surjective. In particular, $\varphi^{-1}(x_1,\ldots,x_{d+1})=\left(\frac{x_1}{1-x_{d+1}}, \frac{x_2}{1-x_{d+1}},\ldots, \frac{x_d}{1-x_{d+1}}, x_{d+1}\right)$ for $x\in \Omega_1$. Hence, $\varphi_1$ is clearly a $C^\infty$--diffeomorphism.

The proof is analogous in the case of $\varphi_2$.
\end{proof}

So, we immediately obtain

\begin{lemma}\label{l.d-isometry} Let $\Phi_1\colon C(\Omega_1)\to C(T_{{d+1,1-\delta}})$  and $\Phi_2\colon C(\Omega_2)\to C(T_{{1-\delta,d+1}})$ be the operators defined by
\[
\Phi_i(u)=u\circ\varphi_i,\quad u\in C(\Omega_i),\, i=1,2.
\]
Then $\Phi_i$ is a surjective isometry for  $i=1,2$. In particular,  $\Phi_1(C^m(T_{{d+1,1-\delta}}))=C^(T_{{d+1,1-\delta}})$ and $\Phi_2(C^m(T_{{1-\delta,d+1}}))=C^m(T_{{1-\delta,d+1}})$ for every $m\in\N$.
\end{lemma}

For each $i\in \{1,2\}$ we define
\begin{equation}\label{e.d-nuovi}
\cA_{d+1,i}:=\Phi_i^{-1}\circ A_{d+1,i}\circ\Phi_i, \quad D(\cA_{d+1,i})=\Phi_i^{-1}(D(A_{d+1,i})),
\end{equation}
with $(A_{d+1,i},D(A_{d+1,i}))$ the second order  differential operators defined in  (\ref{e.d-operatore1}) and (\ref{e.d-operatore2}). Then, by Proposition \ref{p.d-anal} and Lemma \ref{l.d-isometry} the operator $(\cA_{d+1,i}, D(\cA_{d+1,i}))$ generates an analytic semigroup of angle $\pi/2$ on $C(\Omega_i)$  for every $i\in \{1,2\}$. We observe that, for $i=1$ and $v=\Phi_1(u)$ for   some $u\in  \Phi_1^{-1}(D(A_{d+1,1}))$, we have
\begin{eqnarray*}
& &\partial_{r_i}v=(1-r_{d+1})\partial_{x_i}u , \	 \partial_{r_i}^2v=(1-r_{d+1})^2\partial_{x_i}^2u,\quad i=1,\ldots,d\\
& & \partial_{r_ir_j}^2v=(1-r_{d+1})^2\partial_{x_ix_j}^2u,\quad i,j=1,\ldots,d\\
& &\partial_{r_{d+1}}v=-\sum_{i=1}^d r_i\partial_{x_i}u + \partial_{x_{d+1}}u, \\
 & & \partial_{r_{d+1}}^2v=\sum_{i,j=1}^dr_ir_j\partial_{x_ix_j}^2u-2\sum_{i=1}^dr_i\partial_{x_ix_{d+1}}^2u+\partial_{x_{d+1}}^2u.
\end{eqnarray*}
So, we obtain 
\begin{eqnarray*}
& & A_{d+1,1}(\Phi_1(u))(r)= A_{d+1,1} v(r)\\
& &\quad = \frac{1}{2(1-r_{d+1})}\sum_{i,j=1}^dr_i(\delta_{ij}-r_j)(1-r_{d+1})^2\partial_{x_ix_j}^2u + \\
& & \quad +\frac{1}{2}r_{d+1}(1-r_{d+1})\left(\sum_{i,j=1}^dr_ir_j\partial^2_{x_ix_j}u-
 2\sum_{i=1}^dr_i\partial^2_{x_ix_{d+1}}u+\partial_{x_{d+1}}^2u\right)\\
& & =\frac{1}{2}\sum_{i,j=1}^dr_i(1-r_{d+1})[\delta_{ij}-r_j(1-r_{d+1})]\partial_{x_ix_j}^2u+\\
& & -\sum_{i=1}^dr_i(1-r_{d+1})r_{d+1}\partial^2_{x_ix_{d+1}}u+\frac{1}{2}r_{d+1}(1-r_{d+1})\partial_{x_{d+1}}^2u,\quad r\in T_{d+1,1-\delta},
\end{eqnarray*}
and hence, 
\begin{eqnarray*}
\cA_{2,1}u(x)&=&\frac{1}{2}\sum_{i,j=1}^dx_i(\delta_{ij}-x_j)\partial_{x_ix_j}^2u(x)-\sum_{i=1}^dx_ix_{d+1}\partial^2_{x_ix_{d+1}}u(x)+\\
& &\quad\ +\frac{1}{2}x_{d+1}(1-x_{d+1})\partial_{x_{d+1}}^2u(x)\\
&=&\frac{1}{2}\sum_{i,j=1}^{d+1}x_i(\delta_{ij}-x_j)\partial_{x_ix_j}^2u(x),\quad x\in \Omega_1,
\end{eqnarray*}
i.e., $\cA_{d+1}|_{\Omega_1}=\cA_{d+1,1}$. 

On the other hand, for $i=2$ and $v=\Phi_2(u)$ for   some $u\in  \Phi_2^{-1}(D(A_{d+1,2}))$, we have
\begin{eqnarray*}
& & \partial_{r_1}v=\partial_{x_1}u-\sum_{i=2}^{d+1}r_i\partial_{x_i}u,\\
& & \partial_{r_1}^2v=\partial_{x_1}^2u-2\sum_{i=2}^{d+1}r_i\partial_{x_1x_i}^2u+\sum_{i,j=2}^{d+1}r_ir_j\partial_{x_ix_j}^2u,\\
& & \partial_{r_i}v=(1-r_1)\partial_{x_i}u,\ \ \partial_{r_i}^2v=(1-r_1)^2\partial_{x_i}^2u,\quad i=2,\ldots, d+1,\\
& & \partial_{r_ir_j}^2v=(1-r_1)^2\partial_{x_ix_j}^2u,\quad i,j=2,\ldots, d+1,
\end{eqnarray*}
So, we obtain 
\begin{eqnarray*}
& & A_{d+1,2}(\Phi_2(u))(r)= A_{d+1,2} v(r)\\
& &\quad=\frac{1}{2}r_1(1-r_1)\Big(\partial_{x_1}^2u-2\sum_{i=2}^{d+1}r_i\partial_{x_1x_i}^2u+\sum_{i,j=2}^{d+1}r_ir_j\partial_{x_ix_j}^2u\Big) + \\
& & \qquad + \frac{1}{1-r_1}\frac{1}{2}\sum_{i,j=2}^{d+1}r_i(\delta_{ij}-r_j)(1-r_1)^2\partial_{x_ix_j}^2u\\
& &\quad =\frac{1}{2}r_1(1-r_1)\partial_{x_1}^2u-\sum_{i=2}^{d+1}r_1r_i(1-r_1)\partial_{x_1x_i}^2u+\\
& & \quad +\frac{1}{2}\sum_{i,j=2}^{d+1}r_i(1-r_1)[\delta_{ij}-r_j(1-r_1)]\partial_{x_ix_j}^2u,\qquad r\in T_{{1-\delta},d+1},
\end{eqnarray*}
and hence, 
\begin{eqnarray*}
\cA_{d+1,2}u(x)&=&\frac{1}{2}x_1(1-x_1)\partial_{x_1}^2u(x)-\sum_{i=2}^{d+1}x_1x_i\partial_{x_1x_i}^2u(x)+\\
& & +\frac{1}{2}\sum_{i,j=2}^{d+1}x_i(\delta_{ij}-x_j)\partial_{x_ix_j}^2u(x)\\
& =&\frac{1}{2}\sum_{i,j=1}^{d+1}x_i(\delta_{ij}-x_j)\partial_{x_ix_j}^2u(x),\quad x\in \Omega_2,
\end{eqnarray*}
i.e., $\cA_{d+1}|_{\Omega_2}=\cA_{d+1,2}$. 

We now  may prove the main theorem of this section.

\begin{theorem}\label{t.d-analyticity} Suppose that Hypotheses \ref{Hypo} hold. Then
the closure $(\cA_{d+1}, D(\cA_{d+1}))$ of $(\cA_{d+1}, C^2(S_{d+1}))$ generates an  analytic  $C_0$--semigroup $(T(t))_{t\geq 0}$ of angle $\pi/2$ on $C(S_{d+1})$. The semigroup is compact.
\end{theorem}

\begin{proof}Fix $0<\delta<\frac 1 2$. 
Let $\{\psi_i\}_{i=1,2}\su C_c^\infty(\R^{d+1})$ such that $\sum_{i=1}^2(\psi_i)^2=1$ on $S_{d+1}$ and 
\begin{eqnarray*}
& & \mathrm{supp}(\psi_1) \subseteq \{ x\in \R^{d+1}\mid x_{d+1}< 1-\delta\}\\
& &\mathrm{supp}(\psi_2)\subseteq\{ x\in \R^{d+1}\mid x_1<1-\delta\}.
\end{eqnarray*}
For the sake of simplicity, 
we still denote by $\psi_i$ the restriction of $\psi_i$ to $\Omega_i$, for $i=1,2$. 

By Proposition \ref{p.d-anal} the operators $\cA_{d+1,i}$ defined according to  (\ref{e.d-nuovi}) generate analytic semigroups of angle $\pi/2$. So, if $0<\theta<\pi$ is a fixed angle, we can find two positive constant $C$ and $R$ such that, for $|\lambda|\geq R$ with $|\ar \lambda|<\theta$, the resolvents $R(\lambda, \cA_{2,i})$ exist and satisfies
\begin{equation}\label{e.d-risolvente}
\|R(\lambda, \cA_{2,i})\|\leq \frac{C}{|\lambda|},\quad i=1,2.
\end{equation} 
Then, 
for every $\lambda\in\{z\in\C\mid |\ar z|<\theta\}$ with $|\lambda|\geq R$, we can  define the operator $S(\lambda)\colon C(S_{d+1})\to C(S_{d+1})$ via  
\begin{equation}\label{e.d-dresolventt}
S(\lambda)u= \sum_{i=1}^2 \psi_i R(\lambda, \cA_{d+1,i}) (\psi_i u),\quad u\in C(S_{d+1}),
\end{equation}
and hence, 
\[
\|S(\lambda)\|\leq \frac{3C}{|\lambda|}.
\]
We observe that  the previous considerations on the differential operators $\cA_{d+1,i}$ ensure,  for every $i=1,2$ and $u\in C(S_{d+1})$,  that
\begin{equation}\label{e.2-dcompr}
 \cA_{d+1}(\psi_i R(\lambda, \cA_{d+1,i})(\psi_iu))= \cA_{d+1,i}(\psi_iR(\lambda, \cA_{d+1,i})(\psi_i u)), 
 \end{equation}
and, for every $f,g\in D(\cA_{d+1,i})$, that 
\begin{equation}\label{e.2-dcompp} 
\cA_{d+1,i}(fg)=g(\cA_{d+1,i}f)+ f(\cA_{d+1,i}g)+\sum_{i,j=1}^{d+1}x_i(\delta_{ij}-x_j)\partial_{x_i}f\partial_{x_j}g.
\end{equation}
By (\ref{e.2-dcompr}) and (\ref{e.2-dcompp}) we obtain, for every  $\lambda\in \{z\in\C\mid |\ar z|<\theta\}$ with $|\lambda|\geq R$ and $u\in C(S_{d+1})$, that
\begin{eqnarray*}
(\lambda-\cA_{d+1})S(\lambda)u &=& \lambda S(\lambda)(u)- \sum_{i=1}^2 \cA_{d+1}(\psi_i R(\lambda, \cA_{d+1,i})(\psi_iu)) \\
&=&\lambda S(\lambda)(u)-\sum_{i=1}^2 \cA_{d+1,i}(\psi_i(R(\lambda, \cA_{d+1,i})(\psi_i u))\\
&=&\sum_{i=1}^2\psi_i (\lambda-\cA_{d+1,i})R(\lambda, \cA_{d+1,i})(\psi_i u)*\\
&-& \sum_{i=1}^2 \cA_{d+1,i}(\psi_i) R(\lambda, \cA_{d+1,i})(\psi_i u) +\\
&-& \sum_{i=1}^2\sum_{j,h=1}^{d+1}x_j(\delta_{jh}-x_h)\partial_{x_j}(R(\lambda, \cA_{d+1,i})(\psi_i u))\partial_{x_h}\psi_i\\
&=:& (I + B(\lambda) + C(\lambda))(u).
\end{eqnarray*}
Applying (\ref{e.d-risolvente}) we obtain, for every $\lambda\in\{z\in\C\mid |\ar z|<\theta\}$ with $|\lambda|\geq R$ and $u\in C(S_{d+1})$, that
 \[
 \|B(\lambda) u\|_{S_{d+1}}\leq \frac{M}{|\lambda|} \|u\|_{S_{d+1}},
 \]
 for some $M>0$.
In order to estimate $C(\lambda)$ we proceed as  follows.

We first observe, for every $f,g\in D(\cA_{d+1})$, that
\begin{eqnarray*}
&  &\sum_{i,j=1}^{d+1}x_i(\delta_{ij}-x_j)\partial_{x_i}f\partial_{x_j}g=\sum_{i=1}^{d+1}x_i(1-x_i)\partial_{x_i}f\partial_{x_i}g+\\
& &\quad -\sum_{i=1}^{d+1}\sum_{j=1,j\not =i}^{d+1}\Big[x_ix_j\partial_{x_i}f\partial_{x_j}g-x_ix_j\partial_{x_i}f\partial_{x_i}g+x_ix_j\partial_{x_i}f\partial_{x_i}g\Big]\\
& & \quad =\sum_{i=1}^{d+1}x_i(1-x_i)\partial_{x_i}f\partial_{x_i}g-\sum_{i=1}^{d+1}\sum_{j=1,j\not=i}^{d+1}x_ix_j\partial_{x_i}f(\partial_{x_j}g-\partial_{x_i}g)+\\
& & \quad -\sum_{i=1}^{d+1}x_i\partial_{x_i}f\partial_{x_i}g\sum_{j=1,j\not=i}^{d+1}x_j\\
& & \quad =\sum_{i=1}^{d+1}x_i\left(1-\sum_{j=1}^{d+1}x_j\right)\partial_{x_i}f\partial_{x_i}g-\sum_{i=1}^{d+1}\sum_{j=1,j\not=i}^{d+1}x_ix_j\partial_{x_i}f(\partial_{x_j}g-\partial_{x_i}g).
\end{eqnarray*}
Next, let $f\in D(\cA_{d+1,1})$ and $v=\Phi_1(f)$. Then $v\in D(A_{d+1,1})$ and the following holds
\begin{equation}\label{e.d-derivate1}
\partial_{r_i}v=(1-r_{d+1})\partial_{x_i}f,\ i=1,\ldots,d,\quad \partial_{r_{d+1}}v=-\sum_{i=1}^dr_i\partial_{x_i}f+\partial_{x_{d+1}}f,
\end{equation}
and hence,
\begin{equation}\label{e.d-derivate11}
\partial_{x_i}f=\frac{1}{1-r_{d+1}}\partial_{r_i}v,\ i=1,\ldots,d,\quad \partial_{x_{d+1}}f=\partial_{r_{d+1}}v+\sum_{i=1}^d\frac{r_i}{1-r_{d+1}}\partial_{r_i}v.
\end{equation}
So, by Proposition \ref{p.d-stimerisol} we obtain, for every $\lambda\in \{z\in\C\mid |\ar z|<\theta\}$ { with } $|\lambda|> R_1:=\max\{R,l_1,l_2\}$, that
\begin{eqnarray*}
& &\left\|\sum_{i,j=1}^{d+1}x_i(\delta_{ij}-x_j)\partial_{x_i}f\partial_{x_j}\psi_1\right\|_{S_{d+1}}\\\
& &=\left\|\sum_{i=1}^{d+1}x_i\left(1-\sum_{j=1}^{d+1}x_j\right)\partial_{x_i}f\partial_{x_i}\psi_1-\sum_{i=1}^{d+1}\sum_{j=1,j\not=i}^{d+1}x_ix_j\partial_{x_i}f(\partial_{x_j}\psi_1-\partial_{x_i}\psi_1)\right\|_{S_{d+1}} \\
&  & \leq c\left(\sum_{i=1}^{d+1}\left\|x_i\left(1-\sum_{j=1}^{d+1}x_j\right)\partial_{x_i}f\right\|_{\Omega_1}+\sum_{i=1}^{d+1}\sum_{j=1,j\not=i}^{d+1}\|x_ix_j\partial_{x_i}f\|_{\Omega_1}\right)\nonumber\\
& & =c\Big(\sum_{i=1}^{d}\left\|r_i(1-r_{d+1})^2\left(1-\sum_{j=1}^dr_j\right)\frac{1}{1-r_{d+1}}\partial_{r_i}v\right\|_{T_{d+1,1-\delta}}+\\
& & + \left\|r_{d+1}(1-r_{d+1})\left(1-\sum_{j=1}^dr_j\right)\left(\partial_{r_{d+1}}v+\sum_{i=1}^d\frac{r_i}{1-r_{d+1}}\partial_{r_i}v\right)\right\|_{T_{d+1,1-\delta}}+\\
& & +\sum_{i=1}^d\sum_{j=1,j\not=i}^d\left\|r_ir_j(1-r_{d+1})^2\frac{1}{1-r_{d+1}}\partial_{r_i}v\right\|_{T_{d+1,1-\delta}}+\\
& &+ \sum_{i=1}^{d}\left\|r_i(1-r_{d+1}))\frac{1}{1-r_{d+1}}\partial_{r_i}v\right\|_{T_{d+1,1-\delta}}+\\
& & +\sum_{i=1}^d\left\|r_{d+1}r_i(1-r_{d+1})\left(\partial_{r_{d+1}}v+\sum_{i=1}^d\frac{r_i}{1-r_{d+1}}\partial_{r_i}v\right)\right\|_{T_{d+1,1-\delta}}\Big)\\
& & \leq \frac{C_1}{\sqrt{|\lambda|}}\|\lambda v-A_{d+1,1}v\|_{T_{d+1,1-\delta}}=\frac{C_1}{\sqrt{|\lambda|}}\|\lambda f-\cA_{d+1,1}f\|_{\Omega_1}\\
& & \leq \frac{C_1}{\sqrt{|\lambda|}}\|\lambda f-\cA_{d+1}f\|_{S_{d+1}}.
\end{eqnarray*}
If $f=R(\lambda, \cA_{d+1,1})(\psi_1u)$ for some $u\in C(S_{d+1})$, it follows, for every $\lambda\in \{z\in \C\mid |\ar z|<\theta\}$ with $|\lambda|>R_1$, that
\begin{equation}\label{e.d-stimaderivate1}
\left\|\sum_{j,h=1}^{d+1}x_j(\delta_{jh}-x_h)\partial_{x_j}(R(\lambda, \cA_{d+1,1})(\psi_1 u))\partial_{x_h}\psi_i\right\|_{S_{d+1}}\leq \frac{C_1}{\sqrt{|\lambda|}}\|u\|_{S_{d+1}}.
\end{equation}
By the simmetry of change of variables, we obtain in analogous way that there exists $C_2>0$ such that, for every $\lambda\in \{z\in \C\mid |\ar z|<\theta\}$ with $|\lambda|>R_1$ and $u\in C(S_{d+1})$, we have
\begin{equation}\label{e.d-stimaderivate2}
\left\|\sum_{j,h=1}^{d+1}x_j(\delta_{jh}-x_h)\partial_{x_j}(R(\lambda, \cA_{d+1,2})(\psi_2 u))\partial_{x_h}\psi_i\right\|_{S_{d+1}}\leq \frac{C_2}{\sqrt{|\lambda|}}\|u\|_{S_{d+1}}.
\end{equation}
Combining  (\ref{e.d-stimaderivate1}) and (\ref{e.d-stimaderivate2}), we obtain that there exists $K>0$ such that, for every $\lambda\in\{z\in \C\mid |\ar z|<\theta\}$ with $|\lambda|>R_1$ and $u\in C(S_{d+1})$,
\[
\|C(\lambda)u\|_{S_{d+1}} \leq \frac{K}{\sqrt{|\lambda|}}\|u\|_{S_{d+1}}.
\]
If $|\lambda|>>1$, then the operator $B=(\lambda-\cA_{d+1})S(\lambda)$ is invertible in $\cL(C(S_{d+1}))$. Hence, there exists $R(\lambda, \cA_{d+1})=S(\lambda)B^{-1}$ and
\begin{equation}\label{e.d-rappresentazione}
\|R(\lambda, \cA_{d+1})\|=\|S(\lambda)B^{-1}\| \leq \frac{C'}{|\lambda|},
\end{equation}
with $C'$ a positive constant independent of $\lambda$, provided $\lambda-\cA_{d+1}$ is injective and, in particular, for $\lambda>0$ as $\cA_{d+1}$ is dissipative. To conclude that the semigroup is analytic of angle $\pi/2$ it now suffices to repeat the argument already used in  the proof of Proposition \ref{p.2-moperator}.

Since $R(\lambda, \cA_{d+1})=S(\lambda)B^{-1}$ for some $\lambda>0$ and the operator $S(\lambda)$ is compact by Proposition \ref{p.d-anal}, the differential operator $(\cA_{d+1},D(\cA_{d+1}))$ has compact resolvent. Thus, the semigroup is also compact, being analytic and hence, norm continuous.
\end{proof}

Recalling that the  eigenvalues of the operator $\cA_{d+1}$ are given by $\lambda_m=-\frac{m(m-1)}{2}$, $m\geq 1$, \cite[Ch.VIII, p.221]{S3}, and using Theorem \ref{t.d-analyticity} together with \cite[Proposition 5.6]{Met} we obtain the following result.

\begin{theorem}\label{t.d-boundedanal} The semigroup generated by $(\cA_{d+1}, D(\cA_{d+1}))$ is bounded analytic of angle $\pi/2$.
\end{theorem}

Moreover, the differential operator $(\cA_{d+1}, D(\cA_{d+1}))$ satisfies

\begin{prop}\label{p.d-gradiente} Suppose that Hypotheses \ref{Hypo} hold. 

Then the closure $(\cA_{d+1}, D(\cA_{d+1}))$  
of the differential operator $(\cA_{d+1},$ $ C^2(S_{d+1}))$ defined in (\ref{e.operator}) 
satisfies the following properties.

(1) There exist $K_b>0$ and $t_b>0$ such that, for every $0<t<t_b$, $i=1,\ldots, d+1$ and $u\in C(S_{d+1})$, we have
\[
\|\sqrt{x_i(1-x_i)}\partial_{x_i} (T(t)u)\|_{S_{d+1}} \leq \frac{K_b}{\sqrt{t}}\|u\|_{S_{d+1}},\
\]
and such that, for every $t\geq t_b$, $i=1,\ldots, d+1$ and $u\in C(S_{d+1})$, we have
\[
\|\sqrt{x_i(1-x_i)}\partial_{x_i} (T(t)u)\|_{S_{d+1}} \leq K_b\|u\|_{S_{d+1}}.
\]

(2) For each $0<\theta<\pi$ there exists two constants $C>0$ and $l>1$ such that, for every $\lambda\in \{z\in \C\mid |\ar z|<\theta\}$ with $|\lambda|>l$, $i=1,\ldots, d+1$ and $u\in C(S_{d+1})$, we have
\[
\|\sqrt{x_i(1-x_i)}\partial_{x_i}(R(\lambda, \cA_{d+1})u)\|_{S_{d+1}}\leq \frac{C}{\sqrt{|\lambda|}}\|u\|_{S_{d+1}},
\]
\end{prop}

\begin{proof} We first prove property (2). According to the notation in the proof of Theorem \ref{t.d-analyticity}, fixed an angle $0<\theta<\pi$  there exists $l>1$ such that,  for every $\lambda\in \{z\in \C\mid |\ar z|<\theta\}$ with $|\lambda|>l$, we have
\[
R(\lambda, \cA_{d+1})=S(\lambda)B^{-1},
\]
where the operators $S(\lambda)$   are defined in (\ref{e.d-dresolventt}) and  $\|B^{-1}\|\leq 2$.
Fix $i\in \{1,\ldots, d+1\}$. So, we obtain,  for every $\lambda\in \{z\in \C\mid |\ar z|<\theta\}$ with $|\lambda|>l$ and $u\in C(S_{d+1})$,  that
\begin{eqnarray}\label{e.d-rgradiente}
& &\|\sqrt{x_i(1-x_i)}\partial_{x_i}(R(\lambda, \cA_{d+1})u)\|_{S_{d+1}} \leq \nonumber\\
& & \quad \leq  \sum_{j=1}^2 \|\sqrt{x_i(1-x_i)}\partial_{x_i}(\psi_i R(\lambda, \cA_{d+1,j}) (\psi_j B^{-1}u))\|_{\Omega_j}\nonumber\\
& & \quad\leq  \sum_{j=1}^2 \Big(\|\sqrt{x_i(1-x_i)}\partial_{x_i}(\psi_j) R(\lambda, \cA_{d+1,j}) (\psi_j B^{-1}u)\|_{\Omega_j}+\nonumber\\
& & \qquad +\|\sqrt{x_i(1-x_i)}\psi_j \partial_{x_i}(R(\lambda, \cA_{d+1,j}) (\psi_j B^{-1}u))\|_{\Omega_j}\Big)\\
& & \leq  c \sum_{j=1}^2\left(\frac{C}{|\lambda|}\|\psi_j B^{-1}u\|_{\Omega_j}+\|\sqrt{x_i(1-x_i)}\psi_j \partial_{x_i}(R(\lambda, \cA_{d+1,j}) (\psi_j B^{-1}u))\|_{\Omega_j}\right)\nonumber\\
&  & \leq  c \sum_{j=1}^2\left(\frac{C}{\sqrt{|\lambda|}}\|B^{-1}\|\|u\|_{S_{d+1}}+\|\sqrt{x_i(1-x_i)} \partial_x(R(\lambda, \cA_{d+1,j}) (\psi_j B^{-1}u))\|_{\Omega_j}\right)\nonumber
\end{eqnarray}
with $c:=\sup_{j=1}^2\|\psi_j\|_{1,\Omega_i}$. To estimate the second addend on the right in (\ref{e.d-rgradiente}) we proceed as follows.

For $j=1$ set $f= R(\lambda, \cA_{d+1,1}) (\psi_1 B^{-1}u)$ and $v=\Phi_1(f)$. Then, by (\ref{e.d-derivate1}) and (\ref{e.d-derivate11}) we have, for $i=1,\ldots,d$, that
\begin{eqnarray}\label{e.d-adde1}
& &\|\sqrt{x_i(1-x_i)}\partial_{x_i}f\|_{\Omega_1}=\left\|\frac{\sqrt{r_i(1-r_{d+1})[1-r_i(1-r_{d+1})]}}{1-r_{d+1}}\partial_{r_{i}}v\right\|_{T_{d+1,1-\delta}}\nonumber\\
& &\quad\leq \frac{C''_1}{\sqrt{|\lambda|}}\|\lambda v-A_{d+1,1}v\|_{T_{d+1,{1-\delta}}}=\frac{C_1}{\sqrt{|\lambda|}}\|\lambda f-\cA_{d+1,1}f\|_{\Omega_1}\nonumber\\
& &\quad=\frac{C_1}{\sqrt{|\lambda|}}\|\psi_1 B^{-1}u\|_{\Omega_1}\leq \frac{C_1}{\sqrt{|\lambda|}}\|B^{-1}\|\|u\|_{S_{d+1}},
\end{eqnarray} 
and we have, for $i=d+1$, that
\begin{eqnarray}\label{e.d-adde2}
& &\|\sqrt{x_{d+1}(1-x_{d+1})}\partial_{x_{d+1}}f\|_{\Omega_1}=\\
& & =\left\|\sqrt{r_{d+1}(1-r_{d+1})}\left(\partial_{r_{d+1}}v+\sum_{i=1}^d\frac{r_i}{1-r_{d+1}}\partial_{r_i}v\right)\right\|_{T_{d+1,1-\delta}}\nonumber\\
& & \leq\left\|\sqrt{r_{d+1}(1-r_{d+1})}\partial_{r_{d+1}}v\right\|_{T_{d+1,1.\delta}}+\nonumber\\
& &\qquad +\sum_{i=1}^d\left\|\sqrt{r_{d+1}(1-r_{d+1})}\frac{r_i}{1-r_{d+1}}\partial_{r_i}v\right\|_{T_{d+1,1-\delta}}\nonumber\\
& &\leq\frac{C''_1}{\sqrt{|\lambda|}}\|\lambda v-A_{d+1,1}v\|_{T_{d+1,{1-\delta}}}=\frac{C'_1}{\sqrt{|\lambda|}}\|\lambda f-\cA_{d+1,1}f\|_{\Omega_1}\nonumber\\
& &=\frac{C'_1}{\sqrt{|\lambda|}}\|\psi_1 B^{-1}u\|_{\Omega_1}\leq \frac{C'_1}{\sqrt{|\lambda|}}\|B^{-1}\|\|u\|_{S_{d+1}},
\end{eqnarray} 
By the simmetry of the change of variables, one analogously shows that there exists $C_2>0$ such that, for every  $\lambda\in \{z\in \C\mid |\ar z|<\theta\}$ with $|\lambda|>l$, $i=1,\ldots, d+1$ and $u\in C(S_{d+1})$,   we have
\begin{equation}\label{e.d-adde3}
\|\sqrt{x_i(1-x_i)}\partial_{x_i}(R(\lambda, \cA_{d+1,2})(\psi_2B^{-1}u))\|_{\Omega_2} \leq \frac{C_2}{\sqrt{|\lambda|}}\|B^{-1}\|\|u\|_{S_{d+1}}.
\end{equation}
Combining \eqref{e.d-rgradiente},  \eqref{e.d-adde1}, \eqref{e.d-adde2} and \eqref{e.d-adde3}  we obtain, for every  $\lambda\in \{z\in \C\mid |\ar z|<\theta\}$ with $|\lambda|>l$, $i=1,\ldots,d+1$ and $u\in C(S_{d+1})$,  that
\[
\|\sqrt{x_i(1-x_i)}\partial_{x_i}(R(\lambda, \cA_{d+1})u)\|_{S_{d+1}} \leq \frac{C}{\sqrt{|\lambda|}}\|u\|_{S_{d+1}}
\]
for some constant $C>0$ independent of $u$ and $\lambda$.

Property (1) follows as in the proof of Proposition \ref{p.2-mgradiente}(1).
\end{proof}

\subsection{The main results}
Finally, we can state and prove the main results of this paper.

\begin{theorem}\label{t.main} The  closure $(\cA_d, D(\cA_d))$ of the operator $(\cA_d, C^2(S_d))$ generates a bounded analytic $C_0$--semigroup of angle $\pi/2$ on $C(S_d)$ for every $d\geq 1$. The semigroup is compact.
\end{theorem}

\begin{proof} The proof is by induction on the integer $d\geq 1$. The case $d=1$ is given  in \cite{CM,Met}. Suppose that the result holds for $d\geq 2$. Then we can apply Theorem \ref{t.d-analyticity} and conclude that the result holds for $d+1$. Thus, the proof is complete.
\end{proof}

Therefore,  a similar argument as in the proof of Proposition \ref{p.2-moperator} together with Theorem \ref{t.main} allow us to show that the following holds.

\begin{theorem} Let $d\geq 1$ and $m$ be a strictly positive function in $C(S_d)$. Then the operator $(m\cA_d, D(\cA_d))$ generates an analytic $C_0$--semigroup of angle $\pi/2$ on $C(S_d)$. The semigroup is  contractive and compact.
\end{theorem} 

\begin{proof} By Theorem \ref{p.2-generation} and \cite[Theorem]{Do} we can conclude that  $(m\cA_d, D(\cA_d))$ generates a contractive $C_0$--semigroup on $C(S_d)$. 
We claim that the semigroup is analytic of angle $\pi/2$. To show this we can proceed as  in the proof of Proposition \ref{p.2-moperator} and hence, we indicate only the main changes. 

 
For each $n\in\N$ let $I_n^j:=[\frac{j-1}{n},\frac{j+1}{n}]$, $j=1, \dots,n-1$, and let $L=\{1,\ldots, n-1\}^d$. Then, for every $j=(j_1,j_2,\ldots, j_d)\in L$,
we define the set
\[
J^j_n=(I_n^{j_1}\times I_n^{j_2}\times \ldots I^{j_d}_n)\cap S_d.
\]
Set $M_n=\{j\in L\mid J^j_n\not =\emptyset\}$ and fix $V_n^j\in J^j_n$ for all  $j\in M_n$.
 Then, we choose $\phi_n^j\in C^\infty(\R^d)$ for all $j\in M_n$ such that supp\,$(\phi_n^j)\subseteq J_n^j$   and  $\sum_{j\in M_n}(\phi_n^j)^2=1$. 
We observe that, if $v_j\in C(Q)$, for $j\in M_n$,  and $x\in S_d$, then there exists $\ov{j}\in M_n$ such that $x\in J_n^{\ov{j}}$ and hence, 
\[
\sum_{j\in M_n}\phi_n^j(x)v_j(x)= \sum_{j\in M_{n,0}}\phi_n^{j}(x)v_j(x),
\]
where $M_{n,0}=\{j=(\ov{j}_1+h_1,\ov{j}_2+h_2,\ldots,\ov{j}_d+h_d)\mid \forall i\in\{1,\ldots, d\}\ k_i\in \{-1,0,1\}\}$ so that $M_{n,0}$ contains exactly $3^d$ elements. 
Therefore, we have 
\begin{equation}\label{dparun}
\|\sum_{j\in M_n}\phi_n^jv_j\|_{S_d} \leq 3^d\sup_{j\in M_n}\|\phi_n^jv_j\|_{S_d}.
\end{equation}
Since $(\cA_d,D(\cA_d))$ generates a bounded analytic semigroup of angle $\pi/2$ on $C(S_d)$, for each $\lambda\in \C\setminus (-\infty,0)$, $n\in\N$ and $j\in M_n$, we can define
\[
 R_{jn}(\lambda)=\left(\lambda - m\left(V_{n}^j\right) \cA_d\right)^{-1},
 \]
and hence, for a  fixed angle $0<\theta<\pi$, there exists $K>0$ such that, for every $\lambda\in \{z\in\C\mid |\ar z|<\theta\}$, $n\in\N$ and $j\in M_n$, we have 
\begin{equation}\label{dnormres}
\|R_{jn}(\lambda)\|=\left[m\left(V_n^j\right)\right]^{-1}\left\|R\left(\frac{\lambda}{m(V_n^j)}, \cA_d\right)\right\|\leq \frac{K}{|\lambda|}.
\end{equation}
Moreover, if we set $\mu=\lambda\left[m\left(V_n^j\right)\right]^{-1}$, then we have 
\[
\cA_dR_{jn}(\lambda)=
\left[m\left(V_n^j\right)\right]^{-1}(-I+\lambda R_{jn}(\lambda))
\]
and hence, 
\begin{equation}\label{dnormres2}
\|\cA_d R_{jn}(\lambda)\|\leq 
 \left[m\left(V_n^j\right)\right]^{-1}(1+K)\leq \frac{1+K}{m_0}
\end{equation} 
with $m_0=\min_{x\in S_d}m(x)>0$. 
We now consider the approximate resolvents of the operator $m \cA_d$ defined by
\[
 S_n(\lambda)u = \sum_{j\in M_n} \phi_{n}^j \cdot R_{jn}(\lambda) (\phi_{n}^ju), \quad u\in C(S_d).
 \]
Combining (\ref{dnormres}) with (\ref{dparun}), we obtain, for every $\lambda\in \{z\in\C\mid |\ar z|<\theta\}$ and $n\in\N$, that 
\begin{equation}\label{e.dapprox}
\|S_n(\lambda)\|\leq \frac{3^dK}{|\lambda|}.
\end{equation}
Since we have, for every $\phi,\, \eta\in D(\cA_d)$  that
\[ 
 \cA_d(\phi \eta)=\eta \cA_d(\phi)+ \phi \cA_2(\eta) + \sum_{i,j=1}^dx_i(\delta_{ij}-x_j)\partial_{x_i}\phi\partial_{x_j}\eta,
\]  
the operators $S_n(\lambda)$ satisfy, for every $u\in C(S_d)$, 
\begin{eqnarray*} 
& & (\lambda-m \cA_d)S_n(\lambda)u =  u+ \sum_{j\in M_n} \phi_n^j  \left(m\left(V_n^j\right)-m\right) \cA_d(R_{jn}(\lambda) (\phi_n^ju)) +\\
& &  -\sum_{j\in M_n} m \cA_d(\phi_n^j) \cdot R_{jn}(\lambda) (\phi_n^j u)  -\sum_{j\in M_n}\sum_{i,h=1}^d x_i(\delta_{ih}-x_h)\partial_{x_i}( R_{jn}(\lambda) (\phi_n^j u))\partial_{x_h}\phi_n^j \\
& &\quad =:
(I+C_1(\lambda)+C_2(\lambda)+C_3(\lambda))u.
\end{eqnarray*}
We fix   ${\overline{n}}\in\N$ such that $\sup_{J_{\overline{n}}^j}|m(x)-m(V_{\overline{n}}^j)| \leq \varepsilon=:\frac{m_0}{2.3^d(1+K)}$ for $j\in M_{\ov{n}}$. Then,  from  (\ref{dparun}), (\ref{dnormres}) and (\ref{dnormres2}) and Proposition \ref{p.d-gradiente}(2),  and argumenting as in proof of Proposition \ref{p.2-moperator} we obtain, for every $\lambda\in \{z\in\C\mid |\ar z|<\theta|\}$ with $|\lambda|>l$ and $u\in C(S_d)$,  that
\begin{eqnarray*}
& & \|C_1(\lambda)u\|_{S_d}  <\frac 1 2 \|u\|_{S_d},\\
& &\|C_2(\lambda)u\|_{S_d}  \leq \frac{K'}{|\lambda|}\|u\|_{S_d},\\
 & & \|C_3(\lambda)u\|_{S_d} \leq    \frac{K''}{\sqrt{\lambda}}\|u\|_{S_d},
\end{eqnarray*}
for some positive constants $K',\, K''$ independent of $\lambda$ and $u$.  
Now,  if  $|\lambda|$ is large enough, then we get  $\|C_1(\lambda)+C_2(\lambda)+C_3(\lambda)\|<1$ and hence, the operator $B=(\lambda-m A_2)S_{\ov{n}}(\lambda)$ is invertible in $\cL(C(S_d))$.
So, there exists $R(\lambda, m \cA_d)=S_{\ov{n}}(\lambda)B^{-1}$ in $\cL(C(S_d))$  and, by (\ref{e.dapprox}) $\|R(\lambda, m \cA_d)\|= \|S(\lambda)B^{-1}\|\leq \frac{M'}{|\lambda|}$ for some $M'>0$ independent of $\lambda$, provided $\lambda-m\cA_d$ is injective and, in particular, for $\lambda>0$ as $m\cA_d$ is dissipative. To conclude the proof it now suffices to repeat the argument already used in the proof of Proposition \ref{p.2-moperator}.
\end{proof}


\bigskip
\bibliographystyle{plain}

\end{document}